\documentclass[12pt]{article}
\usepackage{amsmath,amssymb,amsfonts,graphicx,wrapfig}

\newcommand{\EHR}{{\mathrm{EHR}}}
\newcommand{\Pois}{{\mathrm{Pois}}}

\newtheorem{theorem}{Theorem}
\newtheorem{lemma}{Lemma}
\newtheorem{corol}{Corollary}

\begin{document}

\begin{center}
{\LARGE On infinite spectra of first order properties of random graphs\\}

\vspace{0.5cm}

{\large M.E.~Zhukovskii\footnote{Moscow Institute of Physics and Technology, department of discrete mathematics.}}
\end{center}

\section{Introduction}
\label{intro}

Asymptotic behavior of first-order properties probabilities of the Erd\H{o}s--R\'{e}nyi random graph $G(n,p)$ have been widely studied in~\cite{Bol}--\cite{complex_graphs}, \cite{Janson}--\cite{Joint}, \cite{Uspekhi}. Let $n\in\mathbb{N}$, $0\leq p\leq 1.$ Consider a set $\Omega_n=\{G=(V_n,E)\}$ of all undirected graphs without loops and multiple edges with a set of vertices $V_n=\{1,2,...,n\}$. {\it Erd\H{o}s--R\'{e}nyi random graph}~\cite{Bol,Janson,Strange,Uspekhi} is a random element $G(n,p)$ on a probability space $(\Omega,\mathcal{F},{\sf P})$ such that it maps $\Omega$ to $\Omega_n$ and its distribution ${\sf P}_{n,p}$ on $\mathcal{F}_n=2^{\Omega_n}$ is defined in the following way:
$$
 {\sf P}_{n,p}(G)=p^{|E|}(1-p)^{C_n^2-|E|}.
$$
Let us denote the event ``$G(n,p)$ follows a property $L$'' by $\{G(n,p)\models L\}$.

The random graph {\it obeys Zero-One Law}, if for any first order property $L$ (see~\cite{Veresh}) the probability ${\sf P}(G(n,p)\models L)$ either tends to $0$ or tends to $1$. In~\cite{Spencer_01}, it was proved that if $p=n^{-\alpha+o(1)}$, $\alpha\in\mathbb{R}_+\setminus\mathbb{Q}$, then $G(n,p)$ obeys Zero-One Law. To avoid trivialities, we shall restrict ourselves to $0<\alpha<1$ (the case $p=O(1/n)$ was studied in~\cite{Spencer_01}). If $\alpha\in\mathbb{Q}\cap(0,1)$, then $G(n,n^{-\alpha})$ does not obey Zero-One Law (see, e.g.,~\cite{Uspekhi}).

In \cite{Zhuk_01_ext}--\cite{Uspekhi}, Zero-One $k$-Law was studied (the random graph obeys {\it Zero-One $k$-Law}, if for any property $L$ which is expressed by a first-order formula with a quantifier depth at most $k$ (see~\cite{Veresh}) the probability ${\sf P}(G(n,p)\models L)$ either tends to $0$ or tends to $1$). Let us remind that a {\it quantifier depth} of a first-order formula is the maximum number of nested quantifiers. We denote a set of all graph properties which are expressed by first order formulae with a quantifier depth at most $k$ by $\mathcal{L}_k$. Moreover, let $\mathcal{L}=\bigcup_{k\in\mathbb{N}}\mathcal{L}_k$ be the set of all first order graph properties.

In 2012, we proved that if $k\geq 3$ and $\alpha\in(0,1/(k-2))$ (see~\cite{Zhuk_01,Zhuk_01_dan}), then $G(n,n^{-\alpha})$ obeys Zero-One $k$-Law. Moreover, in these papers we proved that $G(n,n^{-1/(k-2)})$ does not obey Zero-One $k$-Law. In 2014~\cite{Zhuk_01_ext}, we proved that if $k>3$  and $\alpha=1-\frac{1}{2^{k-1}+\beta}$, $\beta\in(0,\infty)\setminus\mathcal{Q}$, where $\mathcal{Q}$ is the set of all positive fractions with a numerator at most $2^{k-1}$, then $G(n,n^{-{\alpha}})$ obeys Zero-One $k$-Law. Moreover, in the paper it was proved that $G(n,n^{-\alpha})$ does not obey Zero-One $k$-Law, if $\alpha=1-\frac{1}{2^{k-1}+\beta}$, where $\beta\in\{0,1,\ldots,2^{k-1}-2\}$. Finally, in~\cite{Zhuk_max} it was proved that $G(n,n^{-\alpha})$ obeys Zero-One $k$-Law, if $\alpha\in\{1-\frac{1}{2^k-1},1-\frac{1}{2^k}\}$. Thus, $1-\frac{1}{2^k-2}$ --- is the maximum $\alpha$ in $(0,1)$ such that $G(n,n^{-\alpha})$ does not obey Zero-One $k$-Law.

In the presented paper, we prove (see Section~\ref{new}) that in $(1-\frac{1}{2^{k-1}},1)$ there is only a finite number of $\alpha$ such that $G(n,n^{-\alpha})$ does not obey Zero-One $k$-Law.

If the random graph $G(n,n^{-\alpha})$ does not obey Zero-One $k$-Law for some $\alpha\in(0,1)$ and $k\in\mathbb{N}$, then we say that $\alpha$ is in {\it a spectrum of $k$}. Let us remind that in~\cite{Joint} two notions of spectra of a first-order property $L\in\mathcal{L}$ were considered. The first considers
$p=n^{-\alpha}$.  $S^1(L)$ is a set of $\alpha\in (0,1)$ which does {\em not} satisfy the following property: With $p(n)=n^{-\alpha}$, $\lim_{n\rightarrow\infty}{\sf P}(G(n,p(n))\models L)$ exists and is
either zero or one.  The second considers $p=n^{-\alpha+o(1)}$. $S^2(L)$ is a set of $\alpha\in (0,1)$ which does {\em not} satisfy the following property: There exists $\delta\in \{0,1\}$ and
$\epsilon > 0$ so that when $n^{-\alpha-\epsilon} < p(n) <
n^{-\alpha+\epsilon}$, $\lim_{n\rightarrow\infty}
{\sf P}(G(n,p(n))\models L)=\delta$. Let $k\in\mathbb{N}$. Denote unions of $S^1(L)$ and $S^2(L)$ over all $L\in \mathcal{L}_k$ by $S_k^1$ and $S_k^2$ respectively.

In~\cite{Spencer_inf}, it was proved that the sets $S_k^1$ and $S_k^2$ are infinite when $k$ is large enough.
There are, up to tautological equivalence, (see, e.g., \cite{Veresh}) only a finite number of first order sentences with a given quantifier depth. Thus, for $j$ either $1$ or $2$, the set
$S_k^j$ is infinite if and only if there is a single $L$ with quantifier\ depth at most $k$ such that $S^j(L)$ is infinite. Therefore, we always search for one property with an infinite
spectrum when we prove that the spectrum $S_k^j$ is infinite.

It is also known~\cite{Spencer_ord} that all limit points of $S_k^1$ and $S_k^2$ are approached only from above.

In~\cite{Joint}, it was proved that the minimum $k_1$ and $k_2$ such that the sets $S_{k_1}^1$ and $S_{k_2}^2$ are infinite are in the sets $\{4,\ldots,12\}$ and $\{4,\ldots,10\}$ respectively. Moreover, in the same paper we estimate the minimum and the maximum limit points of $S_k^1$, $S_k^2$. Denote sets of limit points of $S_k^1$ and $S_k^2$ by $(S_k^1)'$ and $(S_k^2)'$ respectively. Then
$$
 \min(S_k^1)'\in\left[\frac{1}{k-2},\frac{1}{k-11}\right],\text{ if }k\geq 15,\quad
 \min(S_k^2)'\in\left[\frac{1}{k-1},\frac{1}{k-7}\right],\text{ if }k\geq 10,
$$
$$
 \max(S_k^j)'\in\left[1-\frac{1}{2^{k-13}},1-\frac{1}{2^{k-1}}\right],\quad\quad\quad\text{ if }k\geq 16,\,j\in\{1,2\}.
$$

In the next section, we state new results. We prove them in Section~\ref{proofs}. Some statements on a distribution of small subgraphs in the random graph, which were used in our proofs, are formulated in Section~\ref{small}.

\section{New results}
\label{new}

\begin{theorem}
For any $k\geq 5$, $\frac{1}{\lfloor k/2\rfloor}\in (S_k^1)'$.
\label{new_1}
\end{theorem}

So, we obtain a better upper bound on the minimum limit point of $S_k^1$ for any $k\leq 20$ and a better upper bound on the minimum limit point of $S_k^2$ for all $k\leq 12$. Moreover, Theorem~\ref{new_1} and Zero-One $k$-Law from~\cite{Zhuk_01,Zhuk_01_dan} imply the following statement.

\begin{corol}
The minimum $k$ such that the set $S_k^1$ ($S_k^2$) is infinite equals $4$ or $5$.
\end{corol}

Moreover, we obtain a better lower bound on the maximum limit point of spectra (for small $k$ as well).

\begin{theorem}
For any $k\geq 8$, $1-\frac{1}{2^{k-5}}\in (S_k^1)'$.
\label{new_2}
\end{theorem}

An emptiness of an intersection of $S_k^1$ with $(1-\frac{1}{2^{k-1}},1)$ follows from the result, which is stated below.

\begin{theorem}
Let $k>3$, $b$ be arbitrary natural numbers. Moreover, let $\frac{a}{b}$ be an irreducible positive fraction. Denote $\nu=\max\{1,2^{k-1}-b\}$. Let $a\in\{\nu,\nu+1,\ldots,2^{k-1}\}$,
$\alpha=1-\frac{1}{2^{k-1}+a/b}$. Then the random graph $G(n,n^{-\alpha})$ obeys Zero-One $k$-Law.
\label{new_th_3}
\end{theorem}

\section{Small subgraphs in the random graph}
\label{small}

For an arbitrary graph $G=(E,V)$, set $e(G)=|E|$, $v(G)=|V|$, $\rho(G)=\frac{e(G)}{v(G)}$,
$\rho^{\max}(G)=\max_{H\subseteq G}\rho(H)$ ($\rho(G)$ is called {\it a density of $G$}). Denote the number of copies of $G$ in $G(n,p)$ by $N_G$. Denote the property of containing a copy of $G$ by $L_G$.

\begin{theorem} [\cite{Bol_small, Vince}] If $p=o\left(n^{-1/\rho^{\max}(G)}\right),$
then $\lim\limits_{n\rightarrow\infty}{\sf P}(G(n,p)\models L_G)=0.$
If $n^{-1/\rho^{\max}(G)}=o(p)$, then
$\lim\limits_{n\rightarrow\infty}{\sf P}(G(n,p)\models L_G)=1.$
\label{thr_graphs}
\end{theorem}
In other words, the function $n^{-1/\rho^{\max}(G)}$ is {\it a threshold} (see~\cite{Bol,Janson}) for the property $L_G$.

Let $G$ be {\it a strictly balanced graph} (a density of this graph is greater than a density of any its proper subgraph) with $a(G)$ automorphisms.
\begin{theorem} [\cite{Bol_small}] If $p=n^{-1/\rho(G)}$, then
$N_G\stackrel{d}\longrightarrow\Pois\left(\frac{1}{a(G)}\right).$
\label{pois}
\end{theorem}

Consider arbitrary graphs $G$ and $H$ such that $H\subset G$, $V(H)=\{x_1,...,x_m\}$, $V(G)=\{x_1,...,x_l\}$ and the set $E(G)\setminus (E(H)\cup E(G\setminus H))$ is non-empty. Denote $e(G,H)=e(G)-e(H)$, $v(G,H)=v(G)-v(H)$, $\rho(G,H)=\frac{e(G,H)}{v(G,H)}$, $\rho^{\max}(G,H)=\max_{H\subset K\subseteq G}\rho(K,H)$. Moreover, let $e^{\min}(G,H)$ be the minimum number $e(K,H)$ over all graphs $K$ such that $H\subset K\subseteq G$, $\rho(K,H)=\rho^{\max}(G,H)$ and the set $E(K)\setminus (E(H)\cup E(K\setminus H))$ is non-empty. Consider graphs $\tilde{H},$ $\tilde{G}$, where $V(\tilde{H})=\{\tilde{x}_1,...,\tilde{x}_m\}$, $V(\tilde{G})=\{\tilde{x}_1,...,\tilde{x}_l\}$, $\tilde{H}\subset\tilde{G}$. The graph $\tilde{G}$ is called $(G,(x_1,\ldots,x_m))$-{\it extension of the ordered tuple $(\tilde{x}_1,\ldots,\tilde{x}_m)$}, if
$$
 \{x_{i_1},x_{i_2}\}\in E(G)\setminus E(H) \Rightarrow
 \{\tilde{x}_{i_1},\tilde{x}_{i_2}\}\in E(\tilde{G})\setminus E(\tilde{H}).
$$
The extension is called {\it strict}, if
$$
 \{x_{i_1},x_{i_2}\}\in E(G)\setminus E(H) \Leftrightarrow
 \{\tilde{x}_{i_1},\tilde{x}_{i_2}\}\in E(\tilde{G})\setminus E(\tilde{H}).
$$
Denote the property of containing a $(G,(x_1,\ldots,x_m))$-extension of any ordered tuple of $m$ vertices by $L_{(G,H)}$.

\begin{theorem} [\cite{Spencer_ext}] There exists $0<\varepsilon<K$ such that
$$
 \text{if } p\le\varepsilon n^{-1/\rho^{\max}(G,H)}(\ln n)^{1/e^{\min}(G,H)}, \text{ then }
 \lim_{n\rightarrow\infty}{\sf P}(G(n,p)\models L_{(G,H)})=0;
$$
$$
 \text{if } p\ge K n^{-1/\rho^{\max}(G,H)}(\ln n)^{1/e^{\min}(G,H)}, \text{ then }
 \lim_{n\rightarrow\infty}{\sf P}(G(n,p)\models L_{(G,H)})=1.
$$
\label{thr_extensions}
\end{theorem}
Obviously, for {\it a balanced pair $(G,H)$} (the maximum density $\rho^{\max}(G,H)$ equals $\rho(G,H)$) the quantity $\rho^{\max}(G,H)$ in the statement of Theorem~\ref{thr_extensions} can be replaced by $\rho(G,H)$. In the same way as for graphs, the pair $(G,H)$ is called {\it strictly balanced}, if $\rho(G,H)>\rho(K,H)$ for any graph $K$ such that $H\subset K\subset G$.\\

Fix a number $\alpha\in(0,1)$. Set
$$
 v(G,H)=|V(G)\setminus V(H)|, \,\, e(G,H)=|E(G)\setminus E(H)|,
$$
$$
 f_{\alpha}(G,H)=v(G,H)-\alpha e(G,H).
$$
If for any graph $S$ such that $H\subset S\subseteq G$ the inequality $f_{\alpha}(S,H)>0$ holds, then the pair $(G,H)$ is called \emph{$\alpha$-safe} (see~\cite{Janson,Uspekhi}). %If for any graph $S$ such that $H\subseteq S\subset G$ the inequality $f_{\alpha}(G,S)<0$ holds, then the pair $(G,H)$ is called {\it $\alpha$-rigid} (see~\cite{Uspekhi},\cite{Janson}).
Finally, let us introduce a notion of a maximal pair. Let $\tilde H\subset\tilde G\subset\Gamma$ and $T\subset K$, where $V(T)=\{v_1,\ldots,v_t\}$, $t\leq|V(\tilde G)|.$ The pair $(\tilde G,\tilde H)$ is called \emph{$(K,T)$-maximal in $\Gamma$}, if any ordered tuple $\mathbf{t}$ of $t$ vertices from $V(\tilde G)$ with at least one vertex from $V(\tilde G)\setminus V(\tilde H)$ does not have a strict $(K,(v_1,\ldots,v_t))$-extension $\tilde K$ in $\Gamma$ such that the following properties hold. The intersection of the sets $V(\tilde K)$, $V(\tilde G)$ contains vertices from $\mathbf{t}$ only and any vertex from $V(\tilde K)$ which is not in $\mathbf{t}$ and any vertex from $V(\tilde G)$ which is not in $\mathbf{t}$ are not adjacent. Similarly, the graph $\tilde G$ is called {\it $(K,T)$-maximal in $\Gamma$}, if any ordered tuple $\mathbf{t}$ of $t$ vertices from $V(\tilde G)$ does not have a strict $(K,(v_1,\ldots,v_t))$-extension $\tilde K$ in $\Gamma$ such that the following properties hold. The intersection of the sets $V(\tilde K)$, $V(\tilde G)$ contains vertices from $\mathbf{t}$ only and any vertex from $V(\tilde K)$ which is not in $\mathbf{t}$ and any vertex from $V(\tilde G)$ which is not in $\mathbf{t}$ are not adjacent.

Consider the random graph $G(n,p)$, arbitrary vertices $\tilde x_1,...,\tilde x_m\in V_n$ and a random variable $N^{(K,T)}_{(G,H)}(\tilde x_1,...,\tilde x_m)$ that maps each graph $\mathcal{G}$ from $\Omega_n$ to the number of strict $(G,(x_1,\ldots,x_m))$-extensions $\tilde G$ of $(\tilde x_1,...,\tilde x_m)$ in $\mathcal{G}$ such that the pair $(\tilde G,\tilde G|_{\{\tilde x_1,\ldots,\tilde x_m\}})$ is $(K,T)$-maximal in $\mathcal{G}$ (and $N_{(G,H)}(\tilde x_1,...,\tilde x_m)$ is the number of all $(G,(x_1,\ldots,x_m))$-extensions of $(\tilde x_1,...,\tilde x_m)$ in $\mathcal{G}$). Let us state the result, which was proved in~\cite{Joint}, on an asymptotic behavior of this variable.\\

\begin{theorem} [\cite{Joint}]
Let $0<\alpha_1<\alpha_2<1$. Let a pair $(G,H)$ be $\alpha_2$-safe, $f_{\alpha_1}(K,T)<0$ and $v(T)\leq v(G)$. Let $p\in[n^{-\alpha_2},n^{-\alpha_1}]$. Then a.a.s. for any $\tilde{x}_1,\ldots,\tilde{x}_m$ the inequality $N^{(K,T)}_{(G,H)}(\tilde x_1,...,\tilde x_m)>0$ holds.
\label{maximal_extensions_rigid}
\end{theorem}

If $\tilde H=(\varnothing,\varnothing)$ and $(\tilde G,\tilde H)$ is $(K,T)$-maximal in $\Gamma$, then $\tilde G$ is $(K,T)$-maximal in $\Gamma$. Therefore, we can state a particular case of Theorem~\ref{maximal_extensions_rigid} which considers $(K,T)$-maximal graphs. Let $N^{(K,T)}_{G}$ be a random variable that maps each $\mathcal{G}$ from $\Omega_n$ to the number of $(K,T)$-maximal copies of $G$ in $\mathcal{G}$.

\begin{corol}
Let $0<\alpha_1<\alpha_2<1$. Let $G$ be a strictly balanced graph with $\rho(G)<1/\alpha_2$ and $f_{\alpha_1}(K,T)<0$. If $p\in[n^{-\alpha_2},n^{-\alpha_1}]$, then a.a.s.
the inequality $N^{(K,T)}_{G}>0$ holds.
\label{maximal_graphs_rigid}
\end{corol}

Let us call pairs $(G,(x_1,\ldots,x_m))$ and $(\tilde G,(\tilde x_1,\ldots,\tilde x_m))$, where $\{x_1,\ldots,x_m\}\subset V(G)$ and $\{\tilde x_1,\ldots,$ $\tilde x_m\}\subset\tilde V(G)$, {\it isomorphic}, if the graph $\tilde G$ is a strict $(G,(x_1,\ldots,x_m))$-extension of $(\tilde x_1,\ldots,\tilde x_m)$.\\

Moreover, in our proofs we use a lemma on the existence of a copy of a strictly balanced graph without extensions, which is stated below. A method for obtaining such results is introduced in~\cite{complex_graphs}. Here, we use this method to prove the lemma.\\

Let $H$ be a strictly balanced graph, $(G,H)$ be a strictly balanced pair, $\rho(H)=\rho(G,H)=1/\alpha$. Moreover, let $V(H)=\{h_1,\ldots,h_v\}$, where $v=v(H)$. Let $W$ be a set with the maximum cardinality which contains ordered tuples of $v$ vertices from $V_n$ which satisfy the following property. For any two ordered tuples
$w_1=(x_{i_1},\ldots,x_{i_v}),w_2=(x_{i_{\sigma(1)}},\ldots,x_{i_{\sigma(v)}})\in W$ which coincide as sets, a permutation $\sigma$ of $(h_1,\ldots,h_v)$ does not preserve edges of $H$ (i.e. a mapping $\phi:V(H)\rightarrow V(H)$ such that $\phi(h_i)=h_{\sigma(i)}$, $i\in\{1,\ldots,v\}$, is not an automorphism of $H$). Obviously, $|W|=\frac{n!}{(n-v)!a(H)}$. For any $w\in W$, we denote a set of elements of $w$ by $\overline{w}$. For any $w=(x_{i_1},\ldots,x_{i_v})\in W$, consider an event $A_w$ that some spanning subgraph in $G(n,n^{-\alpha})|_{\overline{w}}$ is isomorphic to $H$ and the corresponding isomorphism maps $x_{i_j}$ to
$h_j$ for each $j\in\{1,\ldots,v\}$.

\begin{lemma}
There exists a subsequence $\{n_i\}_{i\in\mathbb{N}}$ of the sequence of positive integers such that the following property holds. With positive asymptotic probability less than $1$, in $G(n_i,n_i^{-\alpha})$ there exists at least one copy of $H$ and for any $w\in V_{n_i}$ either $\overline{A_w}$ holds or there is no $(G,(h_1,\ldots,h_v))$-extension of $w$ in $G(n_i,n_i^{-\alpha})$.
\label{balanced_graphs}
\end{lemma}

{\it Proof.} Denote $N_H^-(w)=\sum_{\tilde w} I(A_{\tilde w})$, where the summation is taken over all $\tilde w\in W$ which do not intersect $w$. Denote $N_H^+(w)=\sum_{\tilde w} \xi_{\tilde w}$, where the summation is taken over all $\tilde w\in W$ which intersect $w$ such that $\overline{\tilde w}\cap\overline{w}\neq\overline{w}$. The random variable $\xi_{\tilde w}$ is defined in the following way. For any $\mathcal{G}\in\Omega_n$, the equality $\xi_{\tilde w}(\mathcal{G})=1$ holds if and only if $\mathcal{G}$ with edges between any two vertices from $\overline{w}\cap\overline{\tilde w}$ follows $A_w$ (otherwise, $\xi_{\tilde w}(\mathcal{G})=0$).
Set $N_H(w)=N_H^-(w)+N_H^+(w)$.

Denote a probability of the event that in $G(n,n^{-\alpha})$ there exists at least one copy of $H$ and for any ordered tuple $w$ of $v$ vertices from $V_n$ either $A_w$ holds or there is no $(G,(h_1,\ldots,h_v))$-extension of $w$ by $\mu_n$. Then
$$
 {\sf P}(N_H>0)\geq \mu_n={\sf P}(N_G=0)-{\sf P}(N_H=0)\geq {\sf P}(N_H=1,N_G=0).
$$
Theorem~\ref{pois} implies $\lim_{n\rightarrow\infty}{\sf P}(N_H>0)=1-e^{-1/a(H)}$. Finally,
$$
 {\sf P}(N_H=1,N_G=0)=\sum_{w\in W}{\sf P}(N_H=1,N_G=0|A_w){\sf
 P}(A_w)=
$$
$$
 =\sum_{w\in W}{\sf P}(N_H(w)=0,N_{(G,H)}(w)=0|A_w){\sf
 P}(A_w)=
 \sum_{w\in W}{\sf P}(N_H(w)=0,N_{(G,H)}(w)=0){\sf
 P}(A_w)=
$$
$$
 ={\sf P}(N_H(w_0)=0,N_{(G,H)}(w_0)=0)\sum_{w\in W}{\sf
 P}(I_w)\sim\frac{1}{a(H)}{\sf
 P}(N^{-}_H(w_0)=0,N_{(G,H)}(w_0)=0),
$$
where $w_0\in W$ is an arbitrary ordered tuple. Asymptotic equality holds, because Theorem~\ref{thr_graphs} implies that a.a.s. in $G(n,n^{-\alpha})$ there does not exist any subgraph with at most $2v$ vertices and a density greater than $1/\alpha$. The probability ${\sf P}(N^{-}_H(w_0)=0,N_{(G,H)}(w_0)=0)$ converges to some positive number which is less than $1$ (see~\cite{complex_graphs}). Therefore, the lemma is proved.\\

\section{Proofs}
\label{proofs}

First of all, let us introduce some notations.

Let $\mathcal{G}$ be an arbitrary graph. Moreover, let $r,s$ be arbitrary natural numbers. For any vertices $x_1,\ldots,x_s$ of $\mathcal{G}$, we denote a set of all common $r$-neighbors of $x_1,\ldots,x_s$ in $\mathcal{G}$ by $N_r(x_1,\ldots,x_s)$ (we omit $\mathcal{G}$ in this notation when there is no risk of confusion). A {\it $r$-neighbor} of a vertex $x$ is a vertex $y$ such that the minimum length of a path which connects $x$ and $y$ equals $r$ (a {\it length of a path} is a number of edges in it). Set $N(x_1,\ldots,x_s):=N_1(x_1,\ldots,x_s)$.

Moreover, for any two arbitrary vertices $x,y$ of $\mathcal{G}$ and any its subgraphs $A,B$ denote a length of a minimal path in $\mathcal{G}$ which connects $x$ and $y$ by $d_{\mathcal{G}}(x,y)$ ({\it a minimal path} is a path with the minimum length among considered paths). Moreover, we call a path which connects $x$ and some vertex of $A$ {\it a minimal path which connects $x$ and $A$ in $\mathcal{G}$} if its length equals $\min_{y\in B}d_{\mathcal{G}}(x,y)$. Set $d_{\mathcal{G}}(x,A)=d_{\mathcal{G}}(A,x)=\min_{v\in V(A)} d_{\mathcal{G}}(x,v)$, $d_{\mathcal{G}}(A,B)=\min_{v\in V(A)}d_{\mathcal{G}}(v,B)$.

\subsection{Proof of Theorem~\ref{new_1}}

Let $k\geq 5$, $m\in\mathbb{N}$, $\alpha=\frac{1}{\lfloor k/2\rfloor}+\frac{1}{\lfloor k/2\rfloor(m+\lfloor k/2\rfloor-1)}$ and
$p=n^{-\alpha}$.

%\begin{lemma} Существует такое свойство $L$, записываемое с помощью формулы
%первого порядка глубины 5, что ${\sf P}(G(n,p)\models L)$ сходится
%при $n\rightarrow\infty$ и предел отличен от $0$ и $1$.
%\label{property_existence}
%\end{lemma}

%{\it Доказательство.}
Consider a set $\tilde\Omega_n$ of all graphs $\mathcal{G}$ from $\Omega_n$ which follow the properties below.
\begin{enumerate}
\item For any strictly balanced pair $(G,H)$ such that $V(H)=\{h_1,\ldots,h_v\}$, $\rho(G,H)<1/\alpha$,
$v\leq m+\lfloor k/2\rfloor-1$, $v(G)\leq 2(m+\lfloor k/2\rfloor-1)$, any ordered tuple of $v$ vertices has a $(G,(h_1,\ldots,h_v))$-extension in $\mathcal{G}$.

\item %В $\mathcal{G}$ существует копия любого графа $G$,
%количество вершин которого не превосходит $2m+2$, а максимальная
%плотность --- числа $1/\alpha$.
For any $G$ with $v(G)\leq 2(m+\lfloor k/2\rfloor+1)$ and $\rho^{\max}(G)>1/\alpha$, in $\mathcal{G}$ there is no copy of $G$.

\end{enumerate}

Theorem~\ref{thr_graphs} and Theorem~\ref{thr_extensions} imply that ${\sf P}(G(n,p)\in\tilde\Omega_n)\rightarrow 1$ as $n\rightarrow\infty$.

Let $L$ be a first-order property which is expressed by the formula $\exists x_1 \ldots \exists x_{\lfloor k/2\rfloor} \,\,
\varphi(x_1,\ldots,x_{\lfloor k/2\rfloor})$ with the quantifier depth $\max(2\lfloor k/2\rfloor,\lfloor k/2\rfloor+3)\leq k$, where $\varphi(x_1,\ldots,x_{\lfloor k/2\rfloor})=$
$$
 [K(x_1,\ldots,x_{\lfloor k/2\rfloor})\wedge
 (\exists y_1\ldots\exists y_{\lfloor k/2\rfloor}
 \,\, [(y_1\in N(x_1,\ldots,x_{\lfloor k/2\rfloor}))\wedge\ldots\wedge(y_{\lfloor k/2\rfloor}\in N(x_1,\ldots,x_{\lfloor k/2\rfloor}))\wedge
$$ 
$$ 
 K(y_1,\ldots,y_{\lfloor k/2\rfloor})])
 \wedge(\neg(\exists z\,\,[R^{2}_z\wedge\ldots\wedge R^{\lfloor k/2\rfloor}_z
 \wedge(\forall y \,\,((y\in N(x_1,\ldots,x_{\lfloor k/2\rfloor}))\Rightarrow R^{1,2}_z(y)))]
 ))].
$$
Here, we use the following notations:
$$
 K(x_1,\ldots,x_{\lfloor k/2\rfloor})=((x_1\sim x_2)\wedge(x_1\sim x_3)\wedge\ldots\wedge(x_1\sim x_{\lfloor k/2\rfloor})\wedge\ldots\wedge(x_{\lfloor k/2\rfloor-1}\sim x_{\lfloor k/2\rfloor})),
$$
$$
 (y\in N(x_1,\ldots,x_{\lfloor k/2\rfloor}))=((y\sim x_1)\wedge\ldots\wedge(y\sim x_{\lfloor k/2\rfloor})).
$$
For any $1\leq i<j\leq \lfloor k/2\rfloor $,
$$
 R^{i,j}_z(a)=(\exists v\,\,[(v\in N(z,a,x_1,\ldots,x_{i-1},x_{i+1},\ldots,x_{j-1},x_{j+1},\ldots,x_{\lfloor k/2\rfloor }))\wedge(v\nsim x_i)\wedge(v\nsim x_j)]).
$$
For any $2\leq i\leq \lfloor k/2\rfloor $,
$$
 R^{i}_z=(\exists v\,\,[v\in N(z,x_1,\ldots,x_{i-1},x_{i+1},\ldots,x_{\lfloor k/2\rfloor })]).
$$
%Если предикат $R_z(a,b)$ является истинным для вершин $a,b$, то
%будем говорить, что эти вершины являются {\it $z$-смежными}.

Suppose that $\mathcal{G}\in\tilde\Omega_n$ follows $L$. Consider vertices $x_1,\ldots,x_{\lfloor k/2\rfloor }$ such that
$\varphi(x_1,\ldots,x_{\lfloor k/2\rfloor })$ is true. Set $X=\mathcal{G}|_{\{x_1,\ldots,x_{\lfloor k/2\rfloor }\}\cup
N(x_1,\ldots,x_{\lfloor k/2\rfloor })}$, $\chi=|V(X)|-\lfloor k/2\rfloor $, where $N(x_1,\ldots,x_{\lfloor k/2\rfloor })=\{x^1,\ldots,x^{\chi}\}$. Let us prove that $\chi\geq m$. Suppose that $\chi< m$. By the definition of $\tilde\Omega_n$, in $\mathcal{G}$ there are vertices $z,v_1,\ldots,v_{\chi+\lfloor k/2\rfloor -1}$ such that for any $i\in\{1,\ldots,\chi\}$ we have $v_i\in N(x^i,z,x_3,\ldots,x_{\lfloor k/2\rfloor })$, and for any $i\in\{\chi+1,\ldots,\chi+\lfloor k/2\rfloor -1\}$ we have
$v_i\in N(z,x_1,\ldots,x_{i-\chi},x_{i-\chi+2},\ldots,x_{\lfloor k/2\rfloor })$. Indeed, in this case the pair
$(\mathcal{G}|_{\{x_1,\ldots,x_{\lfloor k/2\rfloor },v_1,\ldots,v_{\chi+\lfloor k/2\rfloor -1},z\}\cup
N(x_1,\ldots,x_{\lfloor k/2\rfloor })},X)$ is strictly balanced with the density
$$
 \frac{\lfloor k/2\rfloor (\chi+\lfloor k/2\rfloor -1)}{\chi+\lfloor k/2\rfloor }=\frac{1}{1/\lfloor k/2\rfloor +1/(\lfloor k/2\rfloor (\chi+\lfloor k/2\rfloor -1))}<\frac{1}{\alpha}.
$$
This contradicts the property $L$. Therefore, $\chi\geq m$. Now, let us prove that $\chi=m$. Suppose $\chi>m$. Remove from the set $N(x_1,\ldots,x_{\lfloor k/2\rfloor })$ some vertices in such a way that $m+1$ vertices are in the remainder (but $\lfloor k/2\rfloor$ pairwise adjacent vertices are still in the set). Denote a subgraph in $X$ induced by the union of this remainder with $x_1,\ldots,x_{\lfloor k/2\rfloor }$ by $\tilde X$. Then
$$
 \rho(\tilde
 X)\geq\frac{\lfloor k/2\rfloor (m+1)+\lfloor k/2\rfloor (\lfloor k/2\rfloor -1)}{m+1+\lfloor k/2\rfloor }>1/\alpha.
$$
This contradicts Property 2 in the definition of $\tilde\Omega_n$.

So, $\chi=m$. Let $z$ be a vertex such that the predicate $R_z^{1,2}$ is true for all vertices from $N(x_1,\ldots,x_{\lfloor k/2\rfloor })$, the predicate $R_z^i$ is true for any $i\in\{2,\ldots,\lfloor k/2\rfloor \}$. Then in $\mathcal{G}$ there exist vertices $v_1,\ldots,v_j$ such that $z\in N(v_1,\ldots,v_j)$ and the set $\{x_1,\ldots,x_{\lfloor k/2\rfloor -1}\}\cup N(x_1,\ldots,x_{\lfloor k/2\rfloor })$ can be divided into $j$ subsets $N_1,\ldots,N_j$ in the following way: for any $i\in\{1,\ldots,j\}$ and any vertex $y\in N_i$, $y\sim v_i$ and $v_i$ is adjacent to $\lfloor k/2\rfloor -2$ vertices from $\{x_1,\ldots,x_{\lfloor k/2\rfloor }\}\setminus\{y\}$. Set
$Y=\mathcal{G}|_{\{x_1,\ldots,x_{\lfloor k/2\rfloor },v_1,\ldots,v_j,z\}\cup N(x_1,\ldots,x_{\lfloor k/2\rfloor })}$.
Then
$$
 1/\rho(Y)\leq\frac{\lfloor k/2\rfloor +j+1+m}{\lfloor k/2\rfloor (m+\lfloor k/2\rfloor -1)+m+\lfloor k/2\rfloor -1+j(\lfloor k/2\rfloor -1)}.
$$
Note that the inequality $j<m+\lfloor k/2\rfloor -1$ implies $1/\rho(Y)<\alpha$. Thus, from the definition of $\tilde\Omega_n$ it follows that $j\geq m+\lfloor k/2\rfloor -1$. As $j\leq m+\lfloor k/2\rfloor -1$, the equality $j=m+\lfloor k/2\rfloor -1$ holds. Therefore,
$1/\rho(Y)\leq\frac{2(m+\lfloor k/2\rfloor )}{2\lfloor k/2\rfloor (m+\lfloor k/2\rfloor -1)}=\frac{1}{\lfloor k/2\rfloor }+\frac{1}{\lfloor k/2\rfloor (m+\lfloor k/2\rfloor -1)}=\alpha$,
$1/\rho(X)\leq\frac{m+\lfloor k/2\rfloor }{\lfloor k/2\rfloor (m+\lfloor k/2\rfloor -1)}=\alpha$. Property 2 in the definition of $\tilde\Omega_n$ implies equalities $\rho(X)=\rho(Y)=1/\alpha$. As in $\mathcal{G}$ there is no vertex $z$, which follows the above properties,
the graph $\mathcal{G}$ does not contain a copy $Y$, which, in turn, contains $X$.

In the remaining part of the proof, we will use these notations $X$ and $Y$ for the obtained graphs (the first one is strictly balanced, the second one is balanced, the pair $(Y,X)$ is strictly balanced) with the density $1/\alpha$. Moreover, denote the obtained property of $\mathcal{G}$ (the existence of a copy of $X$ such that no copy of $Y$ contains it) by $\tilde L$. So, we have proved that if $\mathcal{G}\in\tilde\Omega_n$ and $\mathcal{G}$ follows $L$, then $\mathcal{G}$ follows $\tilde L$.

Suppose that $\mathcal{G}\in\tilde\Omega_n$ and $\mathcal{G}$ follows $\tilde L$. Obviously, in this case $\mathcal{G}$ follows $L$ as well.

By Lemma~\ref{balanced_graphs}, there exists a partial limit $\lim_{i\rightarrow\infty}{\sf P}(G(n_i,n_i^{-\alpha})\models \tilde L)=c$, which is not $0$ or $1$. Moreover,
$$
 {\sf P}(G(n_i,n_i^{-\alpha})\models L)\sim {\sf P}(G(n_i,n_i^{-\alpha})\in\tilde\Omega_{n_i},G(n_i,n_i^{-\alpha})\models
 L)=
$$
\begin{equation}
 ={\sf P}(G(n_i,n_i^{-\alpha})\in\tilde\Omega_{n_i},G(n_i,n_i^{-\alpha})\models
 \tilde L)\sim{\sf P}(G(n_i,n_i^{-\alpha})\models \tilde L)=c.
\label{end}
\end{equation}
Since $\frac{1}{\lfloor k/2\rfloor }+\frac{1}{\lfloor k/2\rfloor (m+\lfloor k/2\rfloor -1)}\rightarrow\frac{1}{\lfloor k/2\rfloor }$ as $m\rightarrow\infty$, the theorem is proved.

\subsection{Proof of Theorem~\ref{new_2}}

Let $m\geq 2$ be arbitrary natural numbers, $\alpha=1-\frac{1}{2^{k-5}}+\frac{1}{2^{k-5}m}$ and $p=n^{-\alpha}$.

Consider a set $\tilde\Omega_n$ of all graphs $\mathcal{G}$ from $\Omega_n$ which follow the properties below.
\begin{enumerate}
\item For any strictly balanced pair $(G,H)$ such that $V(H)=\{h_1,\ldots,h_v\}$, $\rho(G,H)<1/\alpha$,
$v\leq (2^{k-5}-1)(m-1)+2$, $v(G)\leq 2(2^{k-5}-1)(m-1)+3$, any ordered tuple of $v$ vertices has a $(G,(h_1,\ldots,h_v))$-extension in $\mathcal{G}$.

\item %В $\mathcal{G}$ существует копия любого графа $G$,
%количество вершин которого не превосходит $2m+2$, а максимальная
%плотность --- числа $1/\alpha$.
For any $G$ with $v(G)\leq 2(2^{k-5}-1)(m+1)+2$ and $\rho^{\max}>1/\alpha$, in $\mathcal{G}$ there is no copy of $G$.

\end{enumerate}

Theorem~\ref{thr_graphs} and Theorem~\ref{thr_extensions} imply that ${\sf P}(G(n,p)\in\tilde\Omega_n)\rightarrow 1$ as $n\rightarrow\infty$.

The property of vertices $x$ and $y$ {\it to be at the distance $i$} (i.e., a length of the minimal path which connects $x$ and $y$ equals $i$) is expressed by the following formula:
$$
 D^*_i(x,y)=D_i(x,y)\wedge\left(\neg\left(\bigvee_{j=1}^{i-1} D_{j}(x,y)\right)\right),
$$
where $D_i(x,y)$ --- is the following formula with the quantifier depth $\lceil\log_2 i\rceil$:
$$
 D_i(x,y)=\exists v \,  (D_{i/2}(x,v)\wedge D_{i/2}(y,v)), \,\text{ if $i$ is even},
$$
$$
 D_i(x,y)=\exists v \, (D_{(i-1)/2}(x,v)\wedge
 D_{(i+1)/2}(y,v)),\,\text{ if $i$ is odd},
$$
and $D_1(x,y)=(x\sim y)$, $D_0(x,y)=(x=y)$. Moreover, set $D_{i,j}^*(x,y,z)=D_i^*(x,z)\wedge D_j^*(z,y)$.

Let $L$ be a first-order property which is expressed by the formula $\exists a\exists b \,\, \varphi(a,b)$ with the quantifier depth $k$, where $\varphi(a,b)=$
$$
 (S(a,b)\wedge[\forall u\, (D_{2^{k-6},2^{k-6}}^*(a,b,u)\Rightarrow R(a,u))]\wedge[\neg(\exists z \,\,((z\neq a)\wedge(\forall u \,(D_{2^{k-6},2^{k-6}}^*(a,b,u)\Rightarrow D_{2^{k-5}}^*(u,z)))))]).
$$
The predicate $S(a,b)=$
$$
(D_{2^{k-5}}^*(a,b)\wedge (\neg(\exists u_1\exists u_2\exists x\,\,
[(u_1\neq u_2)\wedge D_{2^{k-6},2^{k-6}}^*(u_1,u_2,b)\wedge D_{2^{k-6},2^{k-6}}^*(u_1,u_2,a)\wedge
\psi(a,b,u_1,u_2,x)]))),
$$
where
$\psi(a,b,u_1,u_2,x)=$
$$
 \left(\neg\left(\left(\bigvee_{s=2^{k-6}}^{2^{k-5}}\bigvee_{i=1}^{s}\bigvee_{j=2^{k-6}-i}^{2^{k-5}-i}(D_{i,s-i}^*(a,u_1,x)\wedge D_j^*(x,u_2))\right)\right.\vee
 \left.\left(\bigvee_{i=1}^{2^{k-6}}(D_{i,2^{k-6}-i}^*(u_1,b,x)\wedge D_{i}^*(u_2,x))\right)\right)\right)
$$
is true when there do not exist two distinct paths with lengthes at most $2^{k-5}$ which connect the vertex $a$ and two distinct vertices from the set $N_{2^{k-6}}(a,b)$ (moreover, any two distinct vertices from $N_{2^{k-6}}(a,b)$ do not have common neighbors) and there do not exist two distinct intersecting paths with length $2^{k-6}$ which connect the vertex $b$ and two distinct vertices from the set $N_{2^{k-6}}(a,b)$. The truth of the predicate $R(a,u)=$
$$
 \left(\exists x_1\exists x_2\,\,\left[D_{2^{k-6},2^{k-6}}^*(a,u,x_1)\wedge D_{2^{k-7},2^{k-7}}^*(a,u,x_2)\wedge(\neg D_{2^{k-7}}^*(x_1,x_2))\wedge\xi(a,x_1,x_2)\wedge\xi(u,x_1,x_2)\right]\right),
$$
$$
 \xi(a,x_1,x_2)=\left(\neg\left(\exists y\,\left(\bigvee_{i=1}^{2^{k-7}-1}(D_{i,2^{k-6}-i}^*(a,x_1,y)\wedge D_{2^{k-7}-i}^*(y,x_2))\right)\right)\right),%\wedge(\neg(\exists y\,(D_{1,3}^*(u,x_1,y)\wedge D_1^*(y,x_2))))]),
$$
implies the existence of two non-intersecting paths with lengthes $2^{k-6}$ and $2^{k-5}$ which connect the vertex $a$ and $u$.%;
%$$
% R_z(u)=D_8^*(u,z)\wedge\left(\neg\left(\exists u_1\exists u_2\exists x\,\,\left[(u_1\neq u_2)\wedge D_{4,4}^*(u_1,u_2,b)\wedge
% \left(\bigvee_{i=1}^{3}\bigvee_{j=4-i}^{8-i}(D_{i,4-i}^*(a,u_1,x)\wedge D_j^*(x,u_2))\right)\right]\right)\right).
%$$

Suppose that a graph $\mathcal{G}\in\tilde\Omega_n$ follows $L$. Let $a,b$ be vertices such that the formula
$\varphi(a,b)$ is true. Let $X$ be a union of all paths with length $2^{k-5}$ which connect $a$ and $b$ in $\mathcal{G}$. Let $\chi$ be a number of all such paths and $N_{2^{k-6}}(a,b)=\{x^1,\ldots,x^{\chi}\}$. Let us prove that $\chi\geq m$. Suppose that $\chi< m$. By the definition of $\tilde\Omega_n$, in $\mathcal{G}$ there exists a vertex $z$ such that for any $i\in\{1,\ldots,\chi\}$ the property $D_{2^{k-5}}^*(x^i,z)$ holds and there exist $\chi$ paths $P_1,\ldots,P_{\chi}$ with length $2^{k-5}$ connecting $z$ and $x^1,\ldots,x^{\chi}$ respectively such that for any distinct $i,j\in\{1,\ldots,\chi\}$ equality $V(P_i)\cap V(P_j)=\{z\}$ holds. Indeed, if these paths exist, then
the pair $(X\cup P_1\cup\ldots\cup P_{\chi},X)$ is strictly balanced and its density equals
$$
 \frac{2^{k-5}\chi}{(2^{k-5}-1)\chi+1}=\frac{1}{1-1/2^{k-5}+1/(\chi2^{k-5})}<
 \frac{1}{1-1/2^{k-5}+1/(m2^{k-5})}=\frac{1}{\alpha}.
$$
This contradicts the property $L$. Therefore, $\chi\geq m$. Finally, let us prove that $\chi=m$. Suppose that
$\chi>m$. Remove from $X$ paths with length $2^{k-5}$ which connect vertices $a,b$ (without the vertices $a,b$) in such a way that $m+1$ paths remain. Add to the remaining graph paths with length $2^{k-5}$ from $\mathcal{G}$ which connect $a$ and vertices from $N_{2^{k-6}}(a,b)$ (one path for each vertex) such that an intersection of any two of these paths equals $\{a\}$ and an intersection of any of these paths with any path from $X$ contains $a$ and one vertex from $N_{2^{k-6}}(a,b)$ only. Denote the final graph by $\tilde X$. Then $\rho(\tilde X)=\frac{2^{k-4}(m+1)}{2(2^{k-5}-1)(m+1)+2}>1/\alpha$. This contradicts Property
2 in the definition of $\tilde\Omega_n$.

So, $\chi=m$. Let $z\neq a$ be a vertex such that the predicate $D_{2^{k-5}}^*(\cdot,z)$ is true
for all vertices from $N_{2^{k-6}}(a,b)$. Then in $\mathcal{G}$ there exist paths $P_1,\ldots,P_m$ with length $2^{k-5}$ which connect a vertex $z$ with vertices $x^1,\ldots,x^m$ respectively. Suppose that for some $i\in\{1,\ldots,m-1\}$ $P_{i+1}\subseteq P_1\cup\ldots\cup P_i$. Set
$$
 P_{i+1}=(\{x^{i+1},v_1,\ldots,v_{2^{k-5}-1},z\},\{\{x^{i+1},v_1\},\{v_1,v_2\},\ldots,\{v_{2^{k-5}-1},z\}\}).
$$
Then for some $j\in\{1,\ldots,i\}$ the vertex $v_1$ is in $V(P_j)$. Obviously, $v_1\neq x^j$ (otherwise, the predicate $D_{2^{k-5}-1}(x^j,z)$ is true). Suppose that $v_1\nsim x^j$ in $\mathcal{G}$. Then the predicate $D_s(z,v_1)$ is true for some natural $s<2^{k-5}-1$. As $v_1\sim x^{i+1}$, the predicate $D_{s+1}(x^{i+1},z)$ is true as well. This contradicts the truth of the predicate $D_{2^{k-5}}^*(x^{i+1},z)$. Therefore, the vertex $v_1$ is a common neighbor of the vertices $x^{i+1}$ and $x^j$. This contradicts the truth of the predicate $S(a,b)$.  So, for any $i\in\{1,\ldots,m-1\}$ $P_{i+1}\nsubseteq P_1\cup\ldots\cup P_i$. Let us replace the graph $X$ with its union with paths with length $2^{k-5}$ from $\mathcal{G}$ which connect $a$ and vertices from $N_{2^{k-6}}(a,b)$ (one path for each vertex) such that an intersection of any two of these paths equals $\{a\}$ and an intersection of any of these paths with any path from $X$ contains $a$ and one vertex from $N_{2^{k-6}}(a,b)$ only. Consider the sequence of graphs $X_0=X$, $X_1=X\cup P_1$, $X_2=X\cup P_1\cup P_2$, $\ldots$, $X_m=X\cup P_1\cup\ldots\cup P_m$. Set $Y:=X_m$. For any $i\in\{0,\ldots,m-1\}$, the graph $X_{i+1}$ is obtained from the graph $X_i$ by adding $n_i\leq 2^{k-5}-1$ vertices and $e_i\geq n_i+1$ edges. Therefore,
$$
 1/\rho(Y)\leq\frac{2(2^{k-5}-1)m+2+n_1+\ldots+n_m+1}{2^{k-4}m+n_1+\ldots+n_m+m}\leq \alpha,
$$
Equalities hold if and only if $n_i=2^{k-5}-1$ and $e_i=2^{k-5}$ for all $i\in\{0,\ldots,m-1\}$. Therefore, by the definition of the set $\tilde\Omega_n$ these equalities hold and $1/\rho(Y)=1/\rho(X)=\alpha$. As in $\mathcal{G}$ there is no vertex $z$ which follows the above properties, the graph $\mathcal{G}$ does not contain a copy of $Y$, which contains the graph $X$.

As in Theorem~\ref{new_1}, in what follows we exploit the notations $X$ and $Y$ for two obtained graphs with the density $1/\alpha$ (obviously, the graph $X$ and the pair $(Y,X)$ are strictly balanced). Moreover, denote
the obtained property of $\mathcal{G}$ (existence of a copy of $X$ such that any copy of $Y$ does not contain it) by $\tilde L$. We proved that if $\mathcal{G}\in\tilde\Omega_n$ and $\mathcal{G}$ follows $L$, then $\mathcal{G}$ follows $\tilde L$ as well.

Finally, suppose that $\mathcal{G}\in\tilde\Omega_n$ and $\mathcal{G}$ follows $\tilde L$. Then, obviously, $\mathcal{G}$ follows $L$ as well.

By Lemma~\ref{balanced_graphs}, there exists a partial limit
$\lim_{i\rightarrow\infty}{\sf P}(G(n_i,n_i^{-\alpha})\models \tilde L)=c$,
which is not 0 or 1. Moreover, Equation~(\ref{end}) hold.
Since $1-\frac{1}{2^{k-5}}+\frac{1}{2^{k-5}m}\rightarrow1-\frac{1}{2^{k-5}}$ as
$m\rightarrow\infty$, the theorem is proved.

\subsection{Proof of Theorem~\ref{new_th_3}}

We start the proof from the statement of the theorem of Ehrenfeucht in Section~\ref{game}. This theorem is the main tool in proofs of zero-one laws. Then in Section~\ref{constructions} we define some supplementary constructions (cyclic extensions), after which in Section~\ref{properties} we describe asymptotic properties of the random graph which imply the existence of a winning strategy of Duplicator. This strategy is described in Sections~\ref{win_strategy}--\ref{subsub}.

\subsubsection{Ehrenfeucht game}
\label{game}

In this section, we state a particular case of Ehrenfeucht theorem (see~\cite{Ehren}), which holds for graphs. First, let us define Ehrenfeucht game $\EHR(G,H,i)$ on graphs $G,H$ and $i$ rounds (see, e.g.,~\cite{Janson,Uspekhi}). Let $V(G)=\{x_1,...,x_n\},$ $V(H)=\{y_1,...,y_m\}$. In the $\nu\mbox{-}$th round ($1 \leq \nu \leq i$), Spoiler chooses a vertex in any graph (he chooses either $x_{j_{\nu}}\in V(G)$ or
$y_{j'_{\nu}}\in V(H)$). Then Duplicator chooses any vertex in the other graph. If Spoiler chooses in the $\mu\mbox{-}$th round, say, the vertex $x_{j_{\mu}}\in V(G),$ $j_{\mu}=j_{\nu}$ ($\nu<\mu$), then Duplicator must choose the vertex $y_{j'_{\nu}}\in V(H)$. If in this round Spoiler chooses, say, a vertex $x_{j_{\mu}}\in V(G),$ $j_{\mu}\notin\{j_1,...,j_{\mu-1}\}$, then Duplicator must choose a vertex $y_{j'_{\mu}}\in V(H)$ such that $j'_{\mu}\notin\{j'_1,...,j'_{\mu-1}\}$. If he can not do this, Spoiler wins. After the last round vertices $x_{j_1},...,x_{j_{i}}\in V(G)$ and $y_{j'_1},...,y_{j'_{i}}\in V(H)$ are chosen. If some of these vertices coincide, then leave out the copies and consider only distinct vertices: $x_{h_1},...,x_{h_l};$
$y_{h'_1},...,y_{h'_l},$ $l \leq i.$ Duplicator wins if and only if the corresponding subgraphs are isomorphic up to the order of the vertices::
$$
 G|_{\{x_{h_1},...,x_{h_l}\}}\cong
 H|_{\{y_{h'_1},...,y_{h'_l}\}}.
$$

\begin{theorem} [\cite{Ehren}]
For any graphs $G,H$ and any $i\in\mathbb{N}$, Duplicator has a winning strategy in the game $\EHR(G,H,i)$ if and only if for any property $L$ which is expressed by a first-order formula with the quantifier depth at most $i$ either both graphs follow $L$ or both graphs do not follow $L$.
\end{theorem}

It can be easily shown that this theorem has the following corollary related to the zero-one laws (see, e.g., \cite{Uspekhi}).

\begin{theorem}
The random graphs $G(n,p)$ obeys zero-one $k$-law if and only if
$$
 \lim\limits_{n,m\rightarrow\infty}{\sf P}(\mbox{Duplicator has a winning strateg in }\EHR(G(n,p(n)),G(m,p(m)),k))=1.
$$
\label{ehren}
\end{theorem}

\subsubsection{Constructions}
\label{constructions}

Let $m\geq 2$ be an arbitrary natural number. Consider a pair of graphs $(G,H)$ such that $G\supset H$. We say that $G$ is a {\it cyclic $m$-extension of $H$}, if one of the following properties holds.

\begin{itemize}

\item The inequality $m\geq 3$ holds. Moreover, there exists a vertex $x_1$ of $G$ such that
$$
 V(G)\setminus V(H)=\{y_1^1,...,y_{t_1}^1,y_1^2,...,y_{t_2}^2\},
$$
$$
 E(G)\setminus E(H)=\{\{x_1,y_1^1\},\{y_1^1,y_2^1\},...,\{y_{t_1-1}^1,y_{t_1}^1\},
 \{y_{t_1}^1,y_{1}^2\},\{y_1^2,y_2^2\},...,\{y_{t_2-1}^2,y_{t_2}^2\},\{y_{t_2}^2,y_{t_1}^1\}\},
$$
where $t_1+t_2\leq m-1$, $t_1\geq 0$, $t_2\geq 2$ (if $t_1=0$, then the vertex $x_1$ is adjacent to vertices $y^2_1,y^2_{t_2}$). In such a situation, $G$ is the {\it first type} extension. %Вершину $x_1$ в таком случае будем называть {\it расширяющейся}.

\item The inequality $m\geq 2$ holds. Moreover, there exist two distinct vertices $x_1,x_2$ of $G$ such that for some $t\leq m-1$
$$
 G=(V(H)\sqcup\{y_1,...,y_t\},E(H)\sqcup\{\{x_1,y_1\},\{y_1,y_2\},...,\{y_{t-1},y_t\},\{y_t,x_2\}\}).
$$
In such a situation, $G$ is the {\it second type} extension.
%Вершины $x_1,x_2$ в таком случае будем называть {\it расширяющимися}.

\end{itemize}

Let $H\subset G$ be two subgraphs in a graph $\Gamma$. The pair $(G,H)$ is {\it cyclically $m$-maximal in $\Gamma$}, if there are no cyclic $m$-extensions of $G$ in $\Gamma$ which are not cyclic $m$-extensions of $H$.

\subsubsection{Properties which imply the existence of Duplicator's winning strategy}
\label{properties}

Let $k>3$, $b$ be arbitrary natural numbers, $\frac{a}{b}$ be an irreducible positive fraction,
$\alpha=1-\frac{1}{2^{k-1}+a/b}$, $p=n^{-\alpha}$. Moreover, let $a\in\{\max\{1,2^{k-1}-b\},\ldots,2^{k-1}\}$.

Let us define a set of graphs $\mathcal{S}$. A graph $G$ is in $\mathcal{S}$ if and only if it follows three properties below.

\begin{itemize}

\item[{\sf 1)}] In $G$, there are no strictly balanced subgraphs with at most $2^{2k}b$ vertices and a density greater than $1/\alpha$.

\item[{\sf 2)}] Let $\mathcal{H}$ be a set of $\alpha$-safe pairs $(H_1,H_2)$ such that $v(H_1)\leq 2^{2k}b+k2^k$. Let $\mathcal{K}$ be a set of pairs $(K_1,K_2)$ such that $v(K_1)\leq 2^k$, $v(K_2)\leq 2$ and $f_{\alpha}(K_1,K_2)<0$. Then for any pair $(H_1,H_2)\in\mathcal{H}$, $V(H_2)=\{v_1,\ldots,v_h\}$, and for any subgraph $G_2\subset G$, $V(G_2)=\{x_1,\ldots,x_h\}$, in $G$ there exists a strict $(H_1,(v_1,\ldots,v_h))$-extension $G_1$ of the ordered tuple $(x_1,\ldots,x_h)$ such that the pair $(G_1,G_2)$ is $(K_1,K_2)$-maximal in $G$ for any pair $(K_1,K_2)\in\mathcal{K}$.

\item[{\sf 3)}] Let $\mathcal{H}$ be a set of pairs $(H_1,H_2)$ such that $v(H_1)\leq 2^k,$ $v(H_2)\leq 2$ and $f_{\alpha}(H_1,H_2)<0$. Then, for any strictly balanced graph $H$ with at most $2^{2k}b$ vertices and $\rho(H)<1/\alpha$, in $G$ there is a copy of $H$ which is $(H_1,H_2)$-maximal in $G$ for any $(H_1,H_2)\in\mathcal{H}$.

\end{itemize}

By Theorem~\ref{thr_graphs}, Theorem~\ref{maximal_extensions_rigid} and Corollary~\ref{maximal_graphs_rigid}, ${\sf P}(G(n,p)\in\mathcal{S})\rightarrow 1$ as $n\rightarrow\infty$. Therefore, by Theorem~\ref{ehren}, the statement of Theorem~\ref{new_th_3} follows from the existence of a winning strategy of Duplicator in EHR$(G,H,k)$ for all pairs $(G,H)$ such that $G,H\in\mathcal{S}$.

\subsubsection{Winning strategy of Duplicator}
\label{win_strategy}

Let $G,H\in\mathcal{S}$. Let $X_r,Y_r$ be chosen in the $r$-th round graphs by Spoiler and Duplicator respectively. So, the sets $\{X_r,Y_r\}$ and $\{G,H\}$ coincide for all $r\in\{1,\ldots,k\}$. We denote vertices which are chosen in the first $r$ rounds in $X_r$ and $Y_r$ by $x^1_r,\ldots,x^r_r$ and $y^1_r,\ldots,y^r_r$ respectively. Let us describe Duplicator's strategy by induction. The strategy is divided into two parts. We denote the first and second strategy by {\bf\sf S} and {\bf\sf SF} respectively. In the first round, Duplicator always use the strategy {\bf\sf S} and follows this strategy until a round such that chosen subgraphs allow to exploit the strategy {\bf\sf SF}, which was introduced in~\cite{Zhuk_ext_proof} (we do not describe this strategy in the presented paper, because its detailed description can be found in~\cite{Zhuk_max}, Section~4.8).

Before describe the strategies, we introduce one more important notion. Let $r$ rounds are finished, $r\in\{1,\ldots,k\}$. Let $l\in\{1,\ldots,r\}$ and graphs ${\tilde X}_{r}^1,\ldots,{\tilde X}_{r}^l\subset
X_{r}$, ${\tilde Y}_{r}^1,\ldots,{\tilde Y}_{r}^l\subset Y_{r}$ which do not have common vertices
satisfy the following properties.

\begin{itemize}

\item[I] The verices $x_{r}^1,\ldots,x_{r}^{r}$ are elements of the set $V({\tilde X}_{r}^1\cup\ldots\cup{\tilde X}_{r}^l)$, the vertices $y_{r}^1,\ldots,y_{r}^{r}$ are elements of the set $V({\tilde Y}_{r}^1\cup\ldots\cup{\tilde Y}_{r}^l)$.

\item[II] For any distinct $j_1,j_2\in\{1,\ldots,l\}$, the inequalities $d_{X_{r}}({\tilde X}_{r}^{j_1},{\tilde X}_{r}^{j_2})>2^{k-r}$, $d_{Y_{r}}({\tilde Y}_{r}^{j_1},{\tilde
Y}_{r}^{j_2})>2^{k-r}$ hold.

\item[III] For any $j\in\{1,\ldots,l\}$, in the graph $X_{r}$ (in the graph $Y_{r}$) there is no cyclic $2^{k-r}$-extension of the graph ${\tilde X}_{r}^j$ (the graph ${\tilde Y}_{r}^j$).

\item[IV] Cardinalities of the sets $V({\tilde X}_{r}^1\cup\ldots\cup{\tilde X}_{r}^l)$, $V({\tilde
Y}_{r}^1\cup\ldots\cup{\tilde Y}_{r}^l)$ are at most $2^{2k}b+2^{k-1}r$.

\item[V] The graphs ${\tilde X}_{r}^j$ and ${\tilde Y}_{r}^j$ are isomorphic for any $j\in\{1,\ldots,l\}$ and there exists a corresponding isomorphism (one for all these pairs of graphs) which maps the vertices $x_{r}^i$ to the vertices $y_{r}^i$, $i\in\{1,\ldots,r\}$.

\end{itemize}

Two ordered tuples of graphs ${\tilde X}_{r}^1,\ldots,{\tilde X}_{r}^l$ and
${\tilde Y}_{r}^1,\ldots,{\tilde Y}_{r}^l$ which follow the above properties we call {\it $(k,r,l)$-regular
equivalent in $(X_{r},Y_{r})$}. Moreover, we denote an isomorphism from Property V by $\varphi(k,r,l)$ (generally speaking, such an isomorphism is not unique, therefore, we consider an arbitrary isomorphism from Property V).\\

Note that $(k,1,1)$-regular equivalence of ${\tilde X}_{1}^1$ and ${\tilde Y}_{1}^1$ is defined by Properties~I, III, IV and V. Moreover, $(k,k,l)$-regular equivalence of ordered tuples ${\tilde X}_{k}^1,\ldots,{\tilde
X}_{k}^l$ and ${\tilde Y}_{k}^1,\ldots,{\tilde Y}_{k}^l$ is defined by Properties~I, II, IV and V.\\

Two graphs $\tilde X_r^1$ and $\tilde Y_r^1$ are called {\it $(k,r)$-equivalent in $(X_r,Y_r)$}, if for $l=1$ Properties I, IV and V hold and in the graph $X_r$ (the graph $Y_r$) there is no cyclic $2^{k-r}-1$-extension of the graph ${\tilde X}_{r}^1$ (the graph ${\tilde Y}_{r}^1$), there is no second type cyclic $2^{k-r}$-extension of the graph $X_r|_{\{x_r^1,\ldots,x_r^r\}}$ (the graph $Y_r|_{\{y_r^1,\ldots,y_r^r\}}$) and there exists at most one cyclic $2^{k-r}$-extension of the graph $\tilde X_r^1$ (the graph $\tilde Y_r^1$).\\

The main idea of Duplicator's strategy is the following. Duplicator should play in such a way that for some $r\in\{1,\ldots,k-1\}$ and $l\in\{1,\ldots,r\}$ in the graphs $X_r$, $Y_r$ $(k,r,l)$-regular equivalent ordered tuples of subgraphs in $(X_r,Y_r)$ are constructed. In the first round, Duplicator must use the strategy {\bf\sf S$_1$} which is described in the next section. After the $r$-th round, $r\in\{1,\ldots,k-3\}$, if $(k,r,l)$-regular equivalent ordered tuples are not constructed, then, as we show, Duplicator either can find
$(k,r)$-equivalent graphs (and then, in the $r+1$-th round, he must use the strategy  {\bf\sf S$_{r+1}$}, which is described in Section~\ref{S2}) or he must use the strategy {\bf\sf S$_{r+1}^1$}, which is described in Section~\ref{S3}. After the strategy {\bf\sf S$_{r+1}^1$}, Duplicator never turns back to the strategy {\bf\sf S$_{r+j}$}, $j\geq 2$. Strategy {\bf\sf SF} is described in~\cite{Zhuk_max} (Section~4.8) and is used by Duplicator in the $r+1$-th round, $r\geq 2$, if and only if after the $r$-th round for some $l\in\{1,\ldots,r\}$ $(k,r,l)$-regular equivalent ordered tuples of graphs in $(X_r,Y_r)$ are constructed. In~\cite{Zhuk_max}, it is proved that Duplicator wins, when he uses the strategy {\bf\sf SF}.

\subsubsection{Strategy {\bf\sf S$_1$}}
\label{S1}

Consider the first round and two possibilities to choose the first vertex by Spoiler.\\

Let in $X_1$ there is no cyclic $2^{k-1}$-extension of $(\{x_1^1\},\varnothing)$. Then Duplicator chooses a vertex $y_1^1\in V(Y_1)$ which satisfies the following property (such a vertex exists because $Y_1\in\mathcal{S}$ and, therefore, $Y_1$ satisfies {\sf 3)}). There are no cyclic $2^{k-1}$-extensions of $(\{y_1^1\},\varnothing)$ in $Y_1$. Set ${\tilde X}^1_1=(\{x_1^1\},\varnothing)$, ${\tilde Y}^1_1=(\{y_1^1\},\varnothing)$. Property III of $(k,1,1)$-regular equivalence of the graphs ${\tilde X}_1^1$ and ${\tilde Y}_1^1$ in $(X_1,Y_1)$ is already proved. Obviously, Properties I, IV and V hold as well. In this case, in the second round Duplicator exploits the strategy {\bf\sf SF}.\\

Let in $X_1$ there exists at least one cyclic $2^{k-1}$-extension ${\tilde X}_1^1$ of $(\{x_1^1\},\varnothing)$. Let us prove that there exists a sequence of graphs $G_1,G_2,\ldots,G_s$ such that

\begin{itemize}

\item[a)] for any $i\in\{1,\ldots,s-1\}$, the graph $G_{i+1}$ is a cyclic $2^{k-1}$-extension of the graph $G_i$ in $X_1$, $G_1$ is a cyclic $2^{k-1}$-extension of the graph $(\{x_1^1\},\varnothing)$,

\item[b)] either $X_1|_{V(G_s)}=G_s$, or $\rho(X_1|_{V(G_s)})<1/\alpha$,

\item[c)] there are no cyclic $2^{k-1}$-extensions of $G_s$ in $X_1$,

\item[d)] if for some $i\in\{1,\ldots,s-1\}$ the graph $G_{i+1}$ is a cyclic $2^{k-1}$-extension of the graph $G_i$, but it is not a cyclic $2^{k-1}-1$-extension of the graph $G_i$, then there exists $\mu\in\{1,\ldots,s-1\}$ such that the graph $G_{\mu+1}$ is a cyclic $2^{k-1}$-extension of the graph $G_{\mu}$, but it is not a cyclic $2^{k-1}-1$-extension of the graph $G_{\mu}$, while in the graph $X_1\setminus (G_{\mu+1}\setminus G_{\mu})$ there is no cyclic $2^{k-1}$-extensions of the graph $G_{\mu}$.

\end{itemize}

Let us prove the existence of such a sequence.

Obviously, there exists a sequence $G_1\subset G_2\ldots\subset G_i$ with the following properties. First, $G_1$ is a cyclic $2^{k-1}$-extension of the graph $(\{x_1^1\},\varnothing)$, $G_j$ is a cyclic $2^{k-1}$-extension of the graph $G_{j-1}$ for any $j\in\{2,\ldots,i\}$. Second, $j=i$ is the first number (if such a number exists) such that $G_j$ is a cyclic $2^{k-1}$-extension of the graph $G_{j-1}$, but it is not a cyclic $2^{k-1}-1$-extension of the graph $G_{j-1}$ (here, $G_0=(\{x_1^1\},\varnothing)$). If such a number does not exist, then there are no cyclic $2^{k-1}$-extensions of $G_i$ in $X_1$ (obviously, $i$ exists and $i\leq 2^{k-1}b+1$, because a density of $G_i$ is greater than $1/\alpha$, if $i=2^{k-1}b+2$, this contradicts Property {\sf 1)}). In the last situation, the sequence $G_1,\ldots,G_s$ ($s=i$), which satisfies Properties a), c) and d), is already built. Nevertheless, if $G_i$ is not the ``last'' extension, then consider an arbitrary cyclic $2^{k-1}$-extension $\hat G_i$ of $G_{i-1}$ in $X_1\setminus(G_i\setminus G_{i-1})$ (if such an extension exists). Let us add cyclic $2^{k-1}$-extensions $\hat G_{i+1},\hat G_{i+2},\ldots$ of previously constructed graphs one by one in a similar way until there are no cyclic $2^{k-1}$-extensions of the graph $\hat G_{\hat s}$ in $X_1\setminus (G_i\setminus G_{i-1})$. Obviously, the graph $\hat G_{\hat s}\cup G_i$ is a cyclic $2^{k-1}$-extension of the graph $\hat G_{\hat s}$, but it is not its cyclic $2^{k-1}-1$-extension. Moreover, there are no cyclic $2^{k-1}$-extensions of $\hat G_{\hat s}$ in $X_1\setminus((\hat G_{\hat s}\cup G_i)\setminus \hat G_{\hat s})$. So, the first $\hat s+1$ graphs of the sequence are constructed: $G_1,\ldots,G_{i-1},\hat G_i,\ldots,\hat G_{\hat s}, \hat G_{\hat s}\cup G_i$. Let us add cyclic $2^{k-1}$-extensions to the graph $\hat G_{\hat s}\cup G_i$ (each next graph is an extension of the previous one) until there are no cyclic $2^{k-1}$-extensions of the final graph in $X_1$. Obviously, we get the sequence of graphs (we denote it by $G_1,G_2,\ldots,G_s$), which follows Properties a), c) and d) (in addition, the inequality $s\leq 2^{k-1}b+1$ holds, because a density of the graph $G_s$ is greater than $1/\alpha$, if $s=2^{k-1}b+2$, this contradicts Property {\sf 1)}).

Suppose that $e(X_1|_{V(G_s)})>e(G_s)$. Moreover, let $e(X_1|_{V(G_s)})-e(G_s)\geq 2$. Since $s\leq 2^{k-1}b+1$, by Property {\sf 1)} the inequalities $\rho^{\max}(X_1|_{V(G_s)})\leq 1/\alpha<1+\frac{1}{2^{k-1}-1}$ hold. Then
$$
 1+\frac{1}{2^{k-1}-1}>\rho^{\max}(X_1|_{V(G_s)})\geq\frac{2^{k-1}+2}{2^{k-1}}=1+\frac{1}{2^{k-2}}.
$$
This contradicts the inequality $k>3$.  So, $e(X_1|_{V(G_s)})-e(G_s)=1$. For any $i\in\{1,\ldots,s\}$, set $e(G_i)-e(G_{i-1})=e_i\leq 2^{k-1}$, where $G_0=(\{x_1^1\},\varnothing)$. Then
$$
 1/\rho(X_1|_{V(G_s)})=\frac{e_1+\ldots+e_s-s+1}{e_1+\ldots+e_s+1}=
 1-\frac{1}{2^{k-1}+\frac{(e_1-2^{k-1})+\ldots+(e_s-2^{k-1})+1}{s}}.
$$
Therefore, either $1/\rho(X_1|_{V(G_s)})=1-\frac{1}{2^{k-1}+\frac{1}{s}}$, or $1/\rho(X_1|_{V(G_s)})\leq 1-\frac{1}{2^{k-1}}<\alpha$, where the last inequality holds, if at least one of $e_i$, $i\in\{1,\ldots,s\}$, is at most $2^{k-1}-1$. In the last case, we arrive at a contradiction with Property {\sf 1)} of the graph $X_1$, because $s\leq 2^{k-1}b+1$. So, $1/\rho(X_1|_{V(G_s)})=1-\frac{1}{2^{k-1}+\frac{1}{s}}$ and $e_1=\ldots=e_s=2^{k-1}$. If $1/\rho(X_1|_{V(G_s)})>\alpha$, then Property b) holds. Moreover, the inequality $1/\rho(X_1|_{V(G_s)})<\alpha$ contradicts Property {\sf 1)} of the graph $X_1$. Therefore, $1-\frac{1}{2^{k-1}+a/b}=1-\frac{1}{2^{k-1}+1/s}$. As $a/b$ --- the irreducible fraction, $a=1$, $b=s$. The last equalities hold only if $2^{k-1}-b\leq 1$. So, $s\geq 7$. Denote vertices of the additional edge, which exists according to our proposition, by $u,v$. Let $u,v\in V(G_{s-1})$. Moreover, let $u\in V(G_{j_1+1})\setminus V(G_{j_1})$, $v\in V(G_{j_2+1})\setminus V(G_{j_2})$, where $0\leq j_1\leq j_2\leq s-2$, $G_0=(\varnothing,\varnothing)$. Obviously, if the set $V(G_{j_2+1})\setminus V(G_{j_2})$ contains more than one vertex, then there exist graphs $\tilde G_{j_2+1},\ldots,\tilde G_{s+1}$ such that for any $j\in\{j_2,\ldots,s\}$ the graph $\tilde G_{j+1}$ is a cyclic $2^{k-1}$-extension of the graph $\tilde G_j$,  where $\tilde G_{j_2}=G_{j_2}$ and for any $j\in\{j_2+2,\ldots,s+1\}$ the equality $\tilde G_j=X_1|_{V(G_{j-1})}$ holds. If $v(G_{j_2+1},G_{j_2})=1$, then
$$
 1+\frac{1}{2^{k-1}-1}>\rho^{\max}(X_1|_{V(G_{j_2+1})})\geq\frac{2^{k-1}+3}{2^{k-1}+1}=1+\frac{1}{2^{k-2}+1/2}.
$$
This contradicts the inequality $k>3$. Obviously, the sequence $G_1,\ldots,G_{i-1},\tilde G_i,\ldots,\tilde G_{s+1}$ follows Properties a)--d) (here, $\tilde G_{s+1}$ is the cyclic $2^{k-1}$-extension of the graph $\tilde G_s$ from Property d)). Finally, let at least one of the vertices $u,v$ (e.g., $v$) is from the set $V(G_s)\setminus V(G_{s-1})$. If the graph $G_s\setminus (G_{s-1}\setminus G_{s-2})$ is a cyclic $2^{k-1}$-extension of the graph $G_{s-2}$ and $u\notin V(G_{s-1})\setminus V(G_{s-2})$, then set $G_{s-1}:=G_s\setminus (G_{s-1}\setminus G_{s-2})$. So, we get the above situation, which is already considered. If either the graph $G_s\setminus (G_{s-1}\setminus G_{s-2})$ is a cyclic $2^{k-1}$-extension of the graph $G_{s-2}$ and $u\in V(G_{s-1})\setminus V(G_{s-2})$, or $G_s\setminus (G_{s-1}\setminus G_{s-2})$ is not a cyclic $2^{k-1}$-extension of the graph $G_{s-2}$, then there exist graphs $\tilde G_s,\tilde G_{s+1}$ such that the graph $\tilde G_s$ is a cyclic $2^{k-1}$-extension of the graph $G_{s-1}$, the graph $\tilde G_{s+1}$ is a cyclic $2^{k-1}$-extension of the graph $\tilde G_s$, $\tilde G_{s+1}=X_1|_{V(G_s)}$, and there are no cyclic $2^{k-1}$-extension of the graph $G_{s-1}$ in $X_1\setminus(\tilde G_s\setminus G_{s-1})$. Therefore, the sequence $G_1,\ldots,G_{s-1},\tilde G_s,\tilde G_{s+1}$ follows Properties a)--d).

So, let $G_1,G_2,\ldots,G_s$ be a sequence which follows Properties a)--d). Let us prove that the graph $X_1|_{V(G_s)}$ is strictly balanced. Let $\tilde G$ be an arbitrary proper subgraph in $X_1|_{V(G_s)}$. Denote $\tilde G_1=X_1|_{V(G_1)}\cap\tilde G$. If $\tilde G_1\neq X_1|_{V(G_1)}$, then $v(\tilde G\cup X_1|_{V(G_1)},\tilde G)\leq 2^{k-1}-1$, $e(\tilde G\cup X_1|_{V(G_1)},\tilde G)\geq v(\tilde G\cup X_1|_{V(G_1)},\tilde G)+1$. Therefore, a density of the graph $\tilde G\cup X_1|_{V(G_1)}$ is at least
$$
 \frac{e(\tilde G)+v(\tilde G\cup X_1|_{V(G_1)},\tilde G)+1}{v(\tilde G)+v(\tilde G\cup X_1|_{V(G_1)},\tilde G)}>
 \min\left\{\frac{e(\tilde G)}{v(\tilde G)},1+\frac{1}{v(\tilde G\cup X_1|_{V(G_1)},\tilde G)}\right\}=\rho(\tilde G).
$$
In the same way, it can be proved that $\rho(X_1|_{V(G_s)})\geq\rho(\tilde G\cup X_1|_{V(G_{s-1})})\geq\ldots\geq\rho(\tilde G\cup X_1|_{V(G_1)})\geq\rho(\tilde G)$, where at least one of the inequalities is strict, because $\tilde G$ is a proper subgraph in $X_1|_{V(G_s)}$. Therefore, the graph $X_1|_{V(G_s)}$ is strictly balanced.

If $\rho(X_1|_{V(G_s)})<1/\alpha$, then set $\tilde X_1^1=X_1|_{V(G_s)}$. By the definition of the graph $Y_1$, it has a subgraph $\tilde Y_1^1$ isomorphic to $\tilde X_1^1$ such that the following property holds. The graph $\tilde Y_1^1$ is $(K,T)$-maximal for any pair $(K,T)$ such that $v(K)\leq 2^k$, $v(T)\leq 2$ and $f_{\alpha}(K,T)<0$.  Let $\varphi:\tilde X_1^1\rightarrow \tilde Y_1^1$ be an isomorphism. Then Duplicator chooses the vertex $y_1^1:=\varphi(x_1^1)$. By the construction of the graphs $\tilde X_1^1$ and $\tilde Y_1^1$, they do not have cyclic $2^{k-1}$-extensions in $X_1$ and $Y_1$ respectively. Therefore, the graphs $\tilde X_1^1$ and $\tilde Y_1^1$ are $(k,1,1)$-regular equivalent in $(X_1,Y_1)$. In the second round Duplicator exploits the strategy {\bf\sf SF}.

Let $\rho(X_1|_{V(G_s)})=1/\alpha$. Then $\rho(G_s)=1/\alpha$ as well. Set $G_0=(\{x_1^1\},\varnothing)$. For any $i\in\{1,\ldots,s\}$, denote $e_i=e(G_i,G_{i-1})$. Then
$$
 1+\frac{1}{2^{k-1}+a/b-1}=\frac{e_1+\ldots+e_s}{e_1+\ldots+e_s-s+1}=1+\frac{1}{\frac{e_1+\ldots+e_s}{s-1}-1}.
$$
Since $a/b$ is the irreducible fraction, $s\geq b+1$. Obviously, the inequality $a\geq\max\{1,2^{k-1}-b\}$ implies the existence of $\mu\in\{0,\ldots,s-1\}$ such that $G_{\mu+1}$ is not a cyclic $2^{k-1}-1$-extension of the graph $G_{\mu}$. Indeed, otherwise
$$
 \rho(G_s)\geq\frac{(2^{k-1}-1)s}{(2^{k-1}-2)s+1}=1+\frac{1}{2^{k-1}-2+\frac{2^{k-1}-1}{s-1}}\geq 1+\frac{1}{2^{k-1}-2+\frac{2^{k-1}-1}{b}}=
$$
$$
 =1+\frac{1}{2^{k-1}+\frac{2^{k-1}-2b-1}{b}}>1/\alpha.
$$
%Если же $a+b$ --- четное число и $a\geq\max\{1,2^{k-2}-b+1\}$, то предположим, что не существует такого %$i\in\{0,\ldots,s-1\}$, что $G_{i+1}$ не является циклическим $2^{k-1}-1$-расширением графа $G_i$. Тогда
%$$
% (2^{k-1}-1)s\geq e_1+\ldots+e_s=2^{k-1}(s-1)+a\frac{s-1}{b}\geq 2^{k-1}(s-1)+(2^{k-2}-b)\frac{s-1}{b}.
%$$
%Поэтому $\frac{s-1}{b}\leq 2-\frac{1}{2^{k-2}}$, а, следовательно, $\frac{s-1}{b}=1$. Окончательно получаем $e_1+\ldots+e_s=2^{k-1}(s-1)+a$. Если для некоторого $i\in\{1,\ldots,s\}$ выполнено $e_i\leq 2^{k-2}$, то $e_1+\ldots+e_s\leq (s-1)(2^{k-1}-1)+2^{k-2}$. Откуда $a\leq 2^{k-2}-b$, получили противоречие. Если, наконец, все $e_1,\ldots,e_s$ --- нечетные числа, то четность числа $e_1+\ldots+e_s$ совпадает как с четностью числа $s$, так и с четностью числа $a$. Но $s=b+1$, что противоречит четности суммы $a+b$. Поэтому среди чисел $e_1,\ldots,e_s$ есть хотя бы одно четное число (обозначим его $e_{\mu+1}$) из $[2^{k-2}+2,2^{k-1}]$.
Since $G_s$ is strictly balanced, $\rho^{\max}(G_{\mu})<1/\alpha$. As $Y_1\in\mathcal{S}$, in $Y_1$ there exists a subgraph $\tilde Y_1^1$ isomorphic to $\tilde X_1^1:=G_{\mu}$ such that the following property holds. The graph $\tilde Y_1^1$ is $(K,T)$-maximal for any pair $(K,T)$ such that $v(K)\leq 2^k$, $v(T)\leq 2$ and $f_{\alpha}(K,T)<0$.  Let $\varphi:\tilde X_1^1\rightarrow \tilde Y_1^1$ be an isomorphism. Then Duplicator chooses the vertex $y_1^1:=\varphi(x_1^1)$. By the construction of the graphs $\tilde X_1^1$ and $\tilde Y_1^1$, they are $(k,1)$-equivalent in $(X_1,Y_1)$. Therefore, in the second round Duplicator exploits the strategy {\bf\sf S$_2$}.

\subsubsection{Strategy {\bf\sf S$_{r+1}$}}
\label{S2}

Let after the $r$-th round, $r\in\{1,\ldots,k-2\}$, there exist graphs $\tilde X_r^1,\tilde Y_r^1$ which are $(k,r)$-equivalent in $(X_r,Y_r)$. Let $\varphi:\tilde X_r^1\rightarrow\tilde Y_r^1$ be an automorphism.

In the $r+1$-th round, Spoiler chooses a vertex $x_{r+1}^{r+1}$. If $X_{r+1}=X_r$, then set $\tilde X_{r+1}^1=\tilde X_r^1$, $\tilde Y_{r+1}^1=\tilde Y_r^1$. Otherwise, set $\tilde X_{r+1}^1=\tilde Y_r^1$, $\tilde Y_{r+1}^1=\tilde X_r^1$.

Let $x_{r+1}^{r+1}\in V(\tilde X_{r+1}^1)$. Duplicator chooses the vertex $y_{r+1}^{r+1}=\varphi(x_{r+1}^{r+1})$, if $X_{r+1}=X_r$, and the vertex $y_{r+1}^{r+1}=\varphi^{-1}(x_{r+1}^{r+1})$, if $X_{r+1}=Y_r$. As in $X_r,Y_r$ there are no cyclic $2^{k-r-1}$-extensions of the graphs $\tilde X_r^1,\tilde Y_r^1$ respectively (by the definition of the $(k,r)$-equivalence), the graphs $\tilde X_{r+1}^1,\tilde X_{r+1}^1$ are $(k,r+1,1)$-regular equivalent in $(X_{r+1},Y_{r+1})$. Therefore, in the $r+2$-th round Duplicator exploits the strategy {\bf\sf SF}.

Let $x_{r+1}^{r+1}\notin V(\tilde X_{r+1}^1)$. Consider two cases: $r<k-2$ and $r=k-2$.

Let $r<k-2$. If $d_{X_{r+1}}(\tilde X_{r+1}^1,x_{r+1}^{r+1})>2^{k-r-1}$ and in $X_{r+1}$ there are no cyclic $2^{k-r-1}$-extensions of the graph $(\{x_{r+1}^{r+1}\},\varnothing)$, then set $\tilde X_{r+1}^2=(\{x_{r+1}^{r+1}\},\varnothing)$. By Property {\sf 2)} of the graph $Y_{r+1}$, it has a vertex $y_{r+1}^{r+1}$ such that $d_{Y_{r+1}}(\tilde Y_{r+1}^1,y_{r+1}^{r+1})=2^{k-r-1}+1$ and there are no cyclic $2^{k-r-1}$-extensions of $(\{y_{r+1}^{r+1}\},\varnothing)$ in $Y_{r+1}$. Set $\tilde Y_{r+1}^2=(\{y_{r+1}^{r+1}\},\varnothing)$. If there is exactly one cyclic $2^{k-r-1}$-extension of $(\{x_{r+1}^{r+1}\},\varnothing)$, then we denote it by $\tilde X_{r+1}^2$. Let $d_{X_{r+1}}(\tilde X_{r+1}^1,\tilde X_{r+1}^2)>2^{k-r-1}$. By Property {\sf 2)} of the graph $\tilde Y_{r+1}^1$, it has a vertex $y_{r+1}^{r+1}$ and a subgraph $\tilde Y_{r+1}^2$ such that $d_{Y_2}(\tilde Y_{r+1}^1,\tilde Y_{r+1}^2)=2^{k-r-1}+1$, pairs $(\tilde Y_{r+1}^2,(\{y_{r+1}^{r+1}\},\varnothing))$ and $(\tilde X_{r+1}^2,(\{x_{r+1}^{r+1}\},\varnothing))$ are isomorphic, and there are no cyclic $2^{k-r-1}$-extensions of $\tilde Y_{r+1}^2$ in $Y_{r+1}$. The property of $(k,r)$-equivalence of the graphs $\tilde X_r^1,\tilde Y_r^1$ in $(X_r,Y_r)$ implies non-existence of cyclic $2^{k-r-1}$-extensions of $\tilde X_r^1$ and $\tilde Y_r^1$ in $X_r$ and $Y_r$ respectively. Obviously, in all the considered cases the ordered tuples $\tilde X_{r+1}^1,\tilde X_{r+1}^2$ and $\tilde Y_{r+1}^1,\tilde Y_{r+1}^2$ are $(k,r+1,2)$-regular equivalent in $(X_{r+1},Y_{r+1})$. Thus, in the $r+2$-th round Duplicator exploits the strategy {\bf\sf SF}. Let $d_{X_{r+1}}(\tilde X_{r+1}^1,\tilde X_{r+1}^2)\leq 2^{k-r-1}$. The property of $(k,r)$-equivalence of the graphs $\tilde X_r^1,\tilde Y_r^1$ in $(X_r,Y_r)$ implies $d_{X_{r+1}}(\tilde X_{r+1}^1,\tilde X_{r+1}^2)=2^{k-r-1}$ and non-existence of cyclic $2^{k-r-1}-1$-extensions of $(\{x_{r+1}^{r+1}\},\varnothing)$ in $X_{r+1}$. In this case, by Property {\sf 2)} of the graph $Y_{r+1}$ Duplicator is able to choose a vertex $y_{r+1}^{r+1}$ such that the following property holds. There exists an isomorphism $L_X\cup\tilde X_{r+1}^1\rightarrow L_Y\cup\tilde Y_{r+1}^1$ which maps the vertices $x_{r+1}^1,\ldots,x_{r+1}^{r+1}$ to the vertices $y_{r+1}^1,\ldots,y_{r+1}^{r+1}$ respectively, where $L_X$ is a minimal path in $X_{r+1}$ which connects $x_{r+1}^{r+1}$ and $\tilde X_{r+1}^1$, $L_Y$ is a minimal path in $Y_{r+1}$ which connects $y_{r+1}^{r+1}$ and $\tilde Y_{r+1}^1$, and the pair $(L_Y\cup\tilde Y_{r+1}^1,\tilde Y_{r+1}^1)$ is cyclically $2^{k-r-1}$-maximal in $Y_{r+1}$. Set $\tilde X_{r+1}^1:=\tilde X_{r+1}^1\cup L_X$, $\tilde Y_{r+1}^1:=\tilde Y_{r+1}^1\cup L_Y$. Next, Dupllicator exploits the strategy {\bf\sf S$_{r+2}^1$}.
Finally, let us prove the the graph $(\{x_{r+1}^{r+1}\},\varnothing)$ has at most one cyclic $2^{k-r-1}$-extension. Indeed, if two such extensions $A$ and $\tilde A$ exist, then
$$
1/\rho(A\cup \tilde A)\leq\frac{2^{k-r-1}+2^{k-r-1}-1}{2^{k-r-1}+2^{k-r-1}}=1-\frac{1}{2^{k-r}}<\alpha.
$$
This contradicts Property {\sf 1)}, since $v(A\cup\tilde A)\leq 2^{k-r}-1$.

Let $d_{X_{r+1}}(\tilde X_{r+1}^1,x_{r+1}^{r+1})\leq 2^{k-r-1}$. Consider a minimal path $L_X$ in $X_{r+1}$ which connects $x_{r+1}^{r+1}$ and $\tilde X_{r+1}^1$. By Property {\sf 2)} of the graph $Y_{r+1}$, there exists a vertex $y_{r+1}^{r+1}$ such that $d_{Y_{r+1}}(\tilde Y_{r+1}^1,y_{r+1}^{r+1})=d_{X_{r+1}}(\tilde X_{r+1}^1,x_{r+1}^{r+1})$, there exists an isomorphism $L_X\cup\tilde X_{r+1}^1\rightarrow L_Y\cup\tilde Y_{r+1}^1$ which maps the vertices $x_{r+1}^1,\ldots,x_{r+1}^{r+1}$ to the vertices $y_{r+1}^1,\ldots,y_{r+1}^{r+1}$ respectively, and the pair $(L_Y\cup\tilde Y_{r+1}^1,\tilde Y_{r+1}^1)$ is cyclically $2^{k-r-1}$-maximal, where $L_Y$ is a minimal path which connects $y_{r+1}^{r+1}$ and $\tilde Y_{r+1}^1$ in $Y_{r+1}$. Obviously, there are no cyclic $2^{k-r-1}$-extensions of the graph $\tilde Y_{r+1}^1$ in $Y_{r+1}$. Set $\tilde Y_{r+1}^1:=\tilde Y_{r+1}^1\cup L_Y$. If there are no cyclic $2^{k-r-1}$-extensions of $L_X\cup \tilde X_{r+1}^1$ in $X_{r+1}$, then set $\tilde X_{r+1}^1:=L_X\cup\tilde X_{r+1}^1$. Obviously, the graphs $\tilde X_{r+1}^1$ and $\tilde Y_{r+1}^1$ are $(k,r+1,1)$-regular equivalent in $(X_{r+1},Y_{r+1})$. Therefore, in the next round Duplicator exploits the strategy {\bf\sf SF}. If there is a cyclic $2^{k-r-1}$-extension of $L_X\cup \tilde X_{r+1}^1$ in $X_{r+1}$, then $d_{X_{r+1}}(x_{r+1}^{r+1},\tilde X_{r+1}^1)=2^{k-r-1}$ and there are no cyclic $2^{k-r-1}-1$-extensions of $L_X\cup\tilde X_{r+1}^1$ in $X_{r+1}$. In this case, the path $L_X$ could be chosen from a set with at most two paths. If there is one such path, then either a cyclic $2^{k-r-1}$-extension of the graph $L_X\cup\tilde X_{r+1}^1$ is the first type extension, or one of the terminal vertices of $L_X$ does not coincide with each of the vertices $x_{r+1}^1,\ldots,x_{r+1}^{r+1}$. If there are two paths, then consider two cases. First, if a cyclic $2^{k-r}$-extension of the graph $\tilde X_{r+1}^1$ is the first type extension, then we choose an arbitrary path $L_X$ from these two paths. Second, if a $2^{k-r}$-extension of the graph $\tilde X_{r+1}^1$ is the second type extension, then $(k,r)$-equivalence of the graphs $\tilde X_{r+1}^1,\tilde Y_{r+1}^1$ in $(X_{r+1},Y_{r+1})$ imply that at least one path does not contain vertices $x_{r+1}^1,\ldots,x_{r+1}^r$. In this case, $L_X$ is such a path. Obviously, the graphs $\tilde X_{r+1}^1:=L_X\cup\tilde X_{r+1}^1,\tilde Y_{r+1}^1$ are $(k,r+1)$-equivalent in $(X_{r+1},Y_{r+1})$. In the next round, Duplicator exploits the strategy {\bf\sf S$_{r+2}$}.\\

Finally, let $r=k-2.$ If $d_{X_{k-1}}(\tilde X_{k-1}^1,x_{k-1}^{k-1})>2$, then set $\tilde X_{k-1}^2=(\{x_{k-1}^{k-1}\},\varnothing)$. By Property {\sf 2)} of the graph $Y_{k-1}$, it contains a vertex $y_{k-1}^{k-1}$ such that $d_{Y_{k-1}}(\tilde Y_{k-1}^1,y_{k-1}^{k-1})=3$. Set $\tilde Y_{k-1}^2=(\{y_{k-1}^{k-1}\},\varnothing)$. Since the graphs $\tilde X_{k-2}^1,\tilde Y_{k-2}^1$ are $(k,k-2)$-equivalent in $(X_{k-2},Y_{k-2})$, they do not have cyclic $2$-extensions in $X_{k-2}$ and $Y_{k-2}$ respectively. Thus, the ordered tuples $\tilde X_{k-1}^1,\tilde X_{k-1}^2$ and $\tilde Y_{k-1}^1,\tilde Y_{k-1}^2$ are $(k,k-1,2)$-regular equivalent in $(X_{k-1},Y_{k-1})$. Therefore, in the $k$-th round Duplicator exploits the strategy {\bf\sf SF}.

If $d_{X_{k-1}}(\tilde X_{k-1}^1,x_{k-1}^{k-1})\leq 2$, then consider a minimal path $L_X$ in $X_{k-1}$, which connects $x_{k-1}^{k-1}$ and $\tilde X_{k-1}^1$. Moreover, let this path connect $x_{k-1}^{k-1}$ and one of the vertices $x_{k-1}^1,\ldots,x_{k-1}^{k-2}$, if such a path with the minimal length exists. By Property {\sf 2)} of the graph $Y_{k-1}$, it contains a vertex $y_{k-1}^{k-1}$ such that $d_{Y_{k-1}}(\tilde Y_{k-1}^1,y_{k-1}^{k-1})=d_{X_{k-1}}(\tilde X_{k-1}^1,x_{k-1}^{k-1})$, there exists an isomorphism $L_X\cup\tilde X_{k-1}^1\rightarrow L_Y\cup\tilde Y_{k-1}^1$ which maps the vertices $x_{k-1}^1,\ldots,x_{k-1}^{k-1}$ to the vertices $y_{k-1}^1,\ldots,y_{k-1}^{k-1}$ respectively, and the pair $(L_Y\cup\tilde Y_{k-1}^1,\tilde Y_{k-1}^1)$ is cyclically $2$-maximal, where $L_Y$ is a minimal path in $Y_{k-1}$ which connects $y_{k-1}^{k-1}$ and $\tilde Y_{k-1}^1$. Obviously, in the graph $Y_{k-1}$ there are no cyclic $2$-extensions of the graph $\tilde Y_{k-1}^1$. If there are no cyclic $2$-extensions of the graph $L_X\cup \tilde X_{k-1}^1$ in $X_{k-1}$, then the graphs $\tilde X_{k-1}^1\cup L_X$ and $\tilde Y_{k-1}^1\cup L_Y$ are $(k,k-1,1)$-regular equivalent in $(X_{k-1},Y_{k-1})$. In the next round, Duplicator exploits the strategy {\bf\sf SF}. If there is a cyclic $2$-extension of the graph $L_X\cup \tilde X_{k-1}^1$ in $X_{k-1}$, then $d_{X_{k-1}}(x_{k-1}^{k-1},\tilde X_{k-1}^1)=2$. Moreover, by the property of $(k,k-2)$-equivalence of the graphs $\tilde X_{k-2}^1,\tilde Y_{k-2}^1$ in $(X_{k-2},Y_{k-2})$, the only path with length $2$ which does not coincide with $L_X$ and connects the vertex $x_{k-1}^{k-1}$ with some vertex of the graph $\tilde X_{k-1}^1$ satisfies the following property. Its terminal vertex (distinct from $x_{k-1}^{k-1}$) either is not one of the vertices $x_{k-1}^1,\ldots,x_{k-1}^{k-2}$, or equals one of the terminal vertices of $L_X$. Obviously, in the $k$-th round, if Spoiler chooses a vertex from one of the graphs $L_X\cup\tilde X_{k-1}^1,L_Y\cup\tilde Y_{k-1}^1$, then Duplicator wins by choosing the image of $x_k^k$ under an isomorphism of the graphs. If Spoiler chooses a vertex outside these graphs which is adjacent to at most one vertex of $x_{k}^1,\ldots,x_k^{k-1}$, then Duplicator has a winning strategy by Property {\sf 2)} of the graph $Y_k$. Obviously, there exist at most two vertices in $\{x_k^1,\ldots,x_k^{k-1}\}$ which are adjacent to $x_k^k$. Finally, if the vertex $x_k^k$ is adjacent to two vertices of $x_k^1,\ldots,x_k^{k-1}$, then Duplicator chooses the vertex with degree $2$ from either the path $L_X$, or the path $L_Y$, and wins.

\subsubsection{Strategy {\bf\sf S$_{r+1}^1$}}
\label{S3}

Let after the $r$-th round, $r\in\{2,\ldots,k-2\}$, there exist induced subgraphs $\tilde X_r^1$ and $\tilde Y_r^1$ of $X_r$ and $Y_r$ respectively such that the following properties hold. The graph $\tilde Y_r^1$ is cyclically $2^{k-r}$-maximal, $x_r^1,\ldots,x_r^r\in V(\tilde X_r^1)$, $y_r^1,\ldots,y_r^r\in V(\tilde Y_r^1)$, there exists an isomorphism $\varphi:\tilde X_r^1\rightarrow\tilde Y_r^1$ which maps the vertices $x_r^1,\ldots,x_r^r$ to the vertices $y_r^1,\ldots,y_r^r$ respectively. Equalities $\tilde X_r^1=\tilde X_{r-1}^1\cup L_X$, $\tilde Y_r^1=\tilde Y_{r-1}^1\cup L_Y$ hold, where $\tilde X_{r-1}^1,\tilde Y_{r-1}^1$ are graphs which have one common vertex with paths $L_X$ and $L_Y$ respectively, $\varphi|_{\tilde X_{r-1}^1}:\tilde X_{r-1}^1\rightarrow\tilde Y_{r-1}^1$ is an isomorphism, the vertices $x_r^1,\ldots,x_r^{r-1}$ are in $V(\tilde X_{r-1}^1)$, the vertices $x_r^r$ and $y_r^r$ are terminal vertices of paths $L_X$ and $L_Y$ and are not from $V(\tilde X_{r-1}^1)$ and $V(\tilde Y_{r-1}^1)$ respectively. Finally, there exists the only cyclic $2^{k-r}$-extension $C_X\cup \tilde X_r^r$ of the graph $\tilde X_r^r$, where $C_X$ is a path with length $l\in[2^{k-r-1},2^{k-r})$ which connects the vertex $x_r^r$ with some not terminal vertex $x$ of the path $L_X$. Moreover, $l+e(L_X)=2^{k-r+1}$ and $d_{X_r}(x,x_r^r)+l=2^{k-r}$.

In the $r+1$-th round, $r\in\{1,\ldots,k-2\}$, Spoiler chooses a vertex $x_{r+1}^{r+1}$. If $X_{r+1}=X_r$, then set $\tilde X_{r+1}^1=\tilde X_r^1$, $\tilde Y_{r+1}^1=\tilde Y_r^1$. Otherwise, set $\tilde X_{r+1}^1=\tilde Y_r^1$, $\tilde Y_{r+1}^1=\tilde X_r^1$ and rename $L_X:=L_Y$, $L_Y:=L_X$.

Let $x_{r+1}^{r+1}\in V(\tilde X_{r+1}^1)$. Duplicator chooses the vertex $y_{r+1}^{r+1}=\varphi(x_{r+1}^{r+1})$, if $X_{r+1}=X_r$. Duplicator chooses the vertex $y_{r+1}^{r+1}=\varphi^{-1}(x_{r+1}^{r+1})$, if $X_{r+1}=Y_r$. There are no cyclic $2^{k-r-1}$-extensions of the graph $\tilde Y_r^1$ in $Y_r$. There is a cyclic $2^{k-r-1}$-extension of the graph $\tilde X_r^1$ if and only if $l=2^{k-r-1}$ (moreover, the number of such extensions does not exceed one). Suppose that the last equality holds.

Let $x_{r+1}^{r+1}\in\tilde X_{r+1}\setminus L_X$. Set $\tilde X_{r+1}^1=\tilde X_{r-1}^1$, $\tilde Y_{r+1}^1=\tilde Y_{r-1}^1$, if $X_{r+1}=X_r$, and $\tilde X_{r+1}^1=\tilde Y_{r-1}^1$, $\tilde Y_{r+1}^1=\tilde Y_{r-1}^1$, otherwise. Set $\tilde X_{r+1}^2=(\{x_{r+1}^{r+1}\},\varnothing)$, $\tilde Y_{r+1}^2=(\{y_{r+1}^{r+1}\},\varnothing)$. Obviously, $d_{X_{r+1}}(\tilde X_{r+1}^1,\tilde X_{r+1}^2)=d_{Y_{r+1}}(\tilde Y_{r+1}^1,\tilde Y_{r+1}^2)=2^{k-r}+2^{k-r-1}>2^{k-r-1}$ and, moreover, there are no cyclic $2^{k-r-1}$-extensions of the graphs $\tilde X_{r+1}^1,\tilde X_{r+1}^2$ in $X_{r+1}$, there are no cyclic $2^{k-r-1}$-extension of the graphs $\tilde Y_{r+1}^1,\tilde Y_{r+1}^2$ in $Y_{r+1}$.

Let $x_{r+1}^{r+1}\in L_X$.

If $x_{r+1}^{r+1}$ and the terminal vertex of the path $L_X$ from $\tilde X_{r+1}^1$ are at a distance less than
$2^{k-r}$, then denote a minimal path which connects $x_{r+1}^{r+1}$ and the vertex from the intersection of $L_X$ and $\tilde X_{r+1}^1$ by $\tilde L_X$. Rename $\tilde X_{r+1}^1:=\tilde X_{r-1}^1\cup \tilde L_X$, $\tilde Y_{r+1}^1:=\varphi(\tilde X_{r-1}^1\cup \tilde L_X)$, if $X_{r+1}=X_r$, and $\tilde X_{r+1}^1:=\tilde Y_{r-1}^1\cup\tilde L_X$, $\tilde Y_{r+1}^1:=\varphi^{-1}(\tilde Y_{r-1}^1\cup L_X)$, otherwise. Set $\tilde X_{r+1}^2=(\{x_{r+1}^{r+1}\},\varnothing)$, $\tilde Y_{r+1}^2=(\{y_{r+1}^{r+1}\},\varnothing)$. Obviously, $d_{X_{r+1}}(\tilde X_{r+1}^1,\tilde X_{r+1}^2)=d_{Y_{r+1}}(\tilde Y_{r+1}^1,\tilde Y_{r+1}^2)>2^{k-r-1}$. Moreover, there are no cyclic $2^{k-r-1}$-extensions of the graphs $\tilde X_{r+1}^1,\tilde X_{r+1}^2$ in $X_{r+1}$ and no cyclic $2^{k-r-1}$-extensions of the graphs $\tilde Y_{r+1}^1,\tilde Y_{r+1}^2$ in $Y_{r+1}$.

If $x_{r+1}^{r+1}$ and the terminal vertex of the path $L_X$ from $\tilde X_{r+1}^1$ are at the distance $d\geq 2^{k-r}$, then denote a minimal path connecting $x_{r+1}^{r+1}$ and the terminal vertex of the path $L_X$ which is not from $\tilde X_{r+1}^1$ by $\tilde L_X$. Rename $\tilde X_{r+1}^1:=\tilde X_{r-1}^1$, $\tilde Y_{r+1}^1:=Y_{r-1}^1$ and set $\tilde X_{r+1}^2=\tilde L_X$, $\tilde Y_{r+1}^2=\varphi(\tilde L_X)$, if $X_{r+1}=X_r$. Rename $\tilde X_{r+1}^1:=\tilde Y_{r-1}^1$, $\tilde Y_{r+1}^1:=X_{r-1}^1$ and set $\tilde X_{r+1}^2=\tilde L_X$, $\tilde Y_{r+1}^2=\varphi^{-1}(\tilde L_X)$, otherwise. Obviously, $d_{X_{r+1}}(\tilde X_{r+1}^1,\tilde X_{r+1}^2)=d_{Y_{r+1}}(\tilde Y_{r+1}^1,\tilde Y_{r+1}^2)\geq 2^{k-r}>2^{k-r-1}$. Moreover, if $d>2^{k-r}$, then there are no cyclic $2^{k-r-1}$-extensions of the graphs $\tilde X_{r+1}^1,\tilde X_{r+1}^2$ in $X_{r+1}$ and no cyclic $2^{k-r-1}$-extensions of the graphs $\tilde Y_{r+1}^1,\tilde Y_{r+1}^2$ in $Y_{r+1}$.

In all the considered cases, ordered tuples $\tilde X_{r+1}^1,\tilde X_{r+1}^2$ and $\tilde Y_{r+1}^1,\tilde Y_{r+1}^2$ are $(k,r+1,2)$-regular equivalent in $(X_{r+1},Y_{r+1})$. Thus, in the next round Duplicator exploits the strategy {\bf\sf SF}.

If in the last case $d=2^{k-r}$, then in the $r+2$-th round Spoiler chooses a vertex $x_{r+2}^{r+2}$ and, next, Duplicator exploits the strategy which is described in Section~\ref{subsub}.

Finally, let $l>2^{k-r-1}$. Then the graphs $\tilde X_{r+1}^1$ and $\tilde Y_{r+1}^1$ are $(k,r+1,1)$-regular equivalent. Thus, in the $r+2$-th round Duplicator exploits the strategy {\bf\sf SF}.\\

If $x_{r+1}^{r+1}\notin V(\tilde X_{r+1}^1)$ but $x_{r+1}^{r+1}$ is in the (only) cyclic $2^{k-r}$-extension of the graph $\tilde X_{r+1}^1$, then denote a minimal path in $X_{r+1}$ which connects $x_{r+1}^r$ and $x_{r+1}^{r+1}$ by $\tilde X_{r+1}^2$. Rename $\tilde X_{r+1}^1:=\tilde X_{r-1}^1$, $\tilde Y_{r+1}^1:=\tilde Y_{r-1}^1$, if $X_{r+1}=X_r$, and $\tilde X_{r+1}^1:=\tilde Y_{r-1}^1$, $\tilde Y_{r+1}^1:=\tilde X_{r-1}^1$, otherwise. By Property {\sf 2)}, in $Y_{r+1}$ there exists a vertex $y_{r+1}^{r+1}$ and a path $\tilde Y_{r+1}^2$ such that $d_{Y_{r+1}}(y_{r+1}^r,\tilde Y_{r+1}^1)=d_{Y_{r+1}}(\tilde Y_{r+1}^2,\tilde Y_{r+1}^1)>2^{k-r-1}$ and the following properties hold. The graphs $\tilde X_{r+1}^2$ and $\tilde Y_{r+1}^2$ are isomorphic, there exists the respective isomorphism which maps the vertex $x_{r+1}^{r+1}$ to the vertex $y_{r+1}^{r+1}$ and in $Y_{r+1}$ there are no cyclic $2^{k-r-1}$-extensions of the graph $\tilde Y_{r+1}^2$. If $d_{X_{r+1}}(x_{r+1}^r,x_{r+1}^{r+1})<2^{k-r-1}$, then, obviously, in $X_{r+1}$ there are no cyclic $2^{k-r-1}$-extensions of $\tilde X_{r+1}^2$. The ordered tuples $\tilde X_{r+1}^1,\tilde X_{r+1}^2$ and $\tilde Y_{r+1}^1,\tilde Y_{r+1}^2$ are $(k,r+1,2)$-regular equivalent in $(X_{r+1},Y_{r+1})$. Thus, in the $r+2$-th round Duplicator exploits the strategy {\bf\sf SF}. If $d_{X_{r+1}}(x_{r+1}^r,x_{r+1}^{r+1})=2^{k-r-1}$, then in the $r+2$-th round Spoiler chooses a vertex $x_{r+2}^{r+2}$ and, next, Duplicator exploits the strategy which is described in Section~\ref{subsub}.\\

Finally, let $x_{r+1}^{r+1}\notin V(\tilde X_{r+1}^1)$ and $x_{r+1}^{r+1}$ be not from a cyclic $2^{k-r}$-extension of the graph $\tilde X_{r+1}^1$. If $d_{X_{r+1}}(x_{r+1}^{r+1},\tilde X_{r+1}^1)\leq 2^{k-r-1}$, then rename $L_X$ in the following way: $L_X$ is a minimal which connects the vertex $x_{r+1}^{r+1}$ and some vertex of the graph $\tilde X_{r+1}^1$. Obviously, in $X_{r+1}$ there are no cyclic $2^{k-r-1}$-extensions of $\tilde X_{r+1}^1\cup L_X$.  By Property {\sf 2)}, in $Y_{r+1}$ there exists a vertex $y_{r+1}^{r+1}$ such that $d_{Y_{r+1}}(\tilde Y_{r+1}^1,y_{r+1}^{r+1})=d_{X_{r+1}}(\tilde X_{r+1}^1,x_{r+1}^{r+1})$ and the following properties hold. There exists an isomorphism $L_X\cup\tilde X_{r+1}^1\rightarrow L_Y\cup\tilde Y_{r+1}^1$ which maps the vertices $x_{r+1}^1,\ldots,x_{r+1}^{r+1}$ to the vertices $y_{r+1}^1,\ldots,y_{r+1}^{r+1}$ respectively and the pair $(L_Y\cup\tilde Y_{r+1}^1,\tilde Y_{r+1}^1)$ is cyclically $2^{k-r-1}$-maximal, where $L_Y$ is a minimal path which connects the vertex $y_{r+1}^{r+1}$ and the graph $\tilde Y_{r+1}^1$ in $Y_{r+1}$. Obviously, the graphs $\tilde X_{r+1}^1:=\tilde X_{r+1}^1\cup L_X$ and $\tilde Y_{r+1}^1:=\tilde Y_{r+1}^1\cup L_Y$ are $(k,r+1,1)$-regular equivalent in $(X_{r+1},Y_{r+1})$. Therefore, in the $r+2$-th round Duplicator exploits the strategy {\bf\sf SF}. If, finally, $d_{X_{r+1}}(x_{r+1}^{r+1},\tilde X_{r+1}^1)> 2^{k-r-1}$, then denote the only (if it exists) cyclic $2^{k-r-1}$-extension of the graph $(\{x_{r+1}^{r+1}\},\varnothing)$ by $\tilde X_{r+1}^2$ (if there are no such extensions, then set  $\tilde X_{r+1}^2=(\{x_{r+1}^{r+1}\},\varnothing)$). The inequality $d_{X_{r+1}}(\tilde X_{r+1}^2,\tilde X_{r+1}^1)>2^{k-r-1}$ holds. By Property {\sf 2)}, in $Y_{r+1}$ there exist a vertex $y_{r+1}^{r+1}$ and a subgraph $\tilde Y_{r+1}^2$ such that $d_{Y_{r+1}}(\tilde Y_{r+1}^2,\tilde Y_{r+1}^1)=2^{k-r-1}+1$, there exists an isomorphism $\tilde X_{r+1}^2\rightarrow\tilde Y_{r+1}^2$ which maps the vertex $x_{r+1}^{r+1}$ to the vertex $y_{r+1}^{r+1}$ and there are no cyclic $2^{k-r-1}$-extensions of the graph $\tilde Y_{r+1}^2$ in $Y_{r+1}$. Obviously, the ordered tuples $\tilde X_{r+2}^1, \tilde X_{r+2}^2$ and $\tilde Y_{r+2}^1,\tilde Y_{r+2}^2$ are $(k,r+1,2)$-regular equivalent in $(X_{r+1},Y_{r+1})$. Therefore, in the $r+2$-th round Duplicator exploits the strategy {\bf\sf SF}.

\subsubsection{The next round strategy}
\label{subsub}

If $X_{r+2}=X_{r+1}$, then set $\tilde X_{r+2}^1=\tilde X_{r+1}^1$, $\tilde X_{r+2}^2=\tilde X_{r+1}^2$, $\tilde Y_{r+2}^1=\tilde Y_{r+1}^1$, $\tilde Y_{r+2}^2=\tilde Y_{r+1}^2$. Otherwise, set $\tilde X_{r+2}^1=\tilde Y_{r+1}^1$, $\tilde X_{r+2}^2=\tilde Y_{r+1}^2$, $\tilde Y_{r+2}^1=\tilde X_{r+1}^1$, $\tilde Y_{r+2}^2=\tilde X_{r+1}^2$. Denote an isomorphism $\tilde X_{r+2}^1\cup\tilde X_{r+2}^2\rightarrow\tilde Y_{r+2}^1\cup\tilde Y_{r+2}^2$ which maps the vertices $x_{r+2}^1,\ldots,x_{r+2}^{r+1}$ to the vertices $y_{r+2}^1,\ldots,y_{r+2}^{r+1}$ respectively by $\varphi$. If $x_{r+2}^{r+2}\in V(\tilde X_{r+2}^1)$, then Duplicator chooses the vertex $y_{r+2}^{r+2}=\varphi(x_{r+2}^{r+2})$. If $r=k-2$, then Duplicator wins. If $r<k-2$, then, obviously, in $X_{r+2}$ there are no cyclic $2^{k-r-2}$-extensions of the graphs $\tilde X_{r+2}^1,\tilde X_{r+2}^2$, in $Y_{r+2}$ there are no $2^{k-r-2}$-extensions of the graphs $\tilde Y_{r+2}^1,\tilde Y_{r+2}^2$. Thus, the ordered tuples $\tilde X_{r+2}^1,\tilde X_{r+2}^2$ and $\tilde Y_{r+2}^1,\tilde Y_{r+2}^2$ are $(k,r+2,2)$-regular equivalent in $(X_{r+2},Y_{r+2})$. Therefore, in the $r+3$-th round Duplicator exploits the strategy {\bf\sf SF}.

If the vertex $x_{r+2}^{r+2}$ is in the only cyclic $2^{k-r-1}$-extension of the graph $\tilde X_{r+2}^2$, then denote a path with minimal length which satisfies the following properties by $\tilde L_X$. Its terminal vertices coincide with the terminal vertices of the path $\tilde X_{r+2}^2$ and the vertex $x_{r+2}^{r+2}$ is in $V(\tilde L_X)$. Obviously, there exists an isomorphism $\tilde\varphi:\tilde X_{r+2}^1\cup L_X\rightarrow Y_{r+2}^1\cup Y_{r+2}^2$ which maps the vertices $x_{r+2}^1,\ldots,x_{r+2}^{r+2}$ to the vertices $y_{r+2}^1,\ldots,y_{r+2}^{r+2}$. Therefore, if $r=k-2$, then Duplicator wins. If $r<k-2$, then, obviously, the ordered tuples $\tilde X_{r+2}^1,\tilde X_{r+2}^2:=\tilde L_X$ and $\tilde Y_{r+2}^1,\tilde Y_{r+2}^2$ are $(k,r+2,2)$-regular equivalent in $(X_{r+2},Y_{r+2})$. Therefore, in the $r+3$-th round Duplicator exploits the strategy {\bf\sf SF}.

If the vertex $x_{r+2}^{r+2}$ is not in the cyclic $2^{k-r-1}$-extension of the graph $\tilde X_{r+2}^2$ and $d_{X_{r+2}}(x_{r+2}^{r+2},\tilde X_{r+2}^1\cup\tilde X_{r+2}^2)\leq 2^{k-r-2}$, then denote a minimal path which connects $x_{r+2}^{r+2}$ and $\tilde X_{r+2}^1\cup\tilde X_{r+2}^2$ by $\tilde L_X$. Obviously, in $X_{r+2}$ there are no cyclic $2^{k-r-2}$-extensions of the graph $\tilde X_{r+2}^1\cup \tilde X_{r+2}^2\cup \tilde L_X$.  By property {\sf 2)}, in $Y_{r+2}$ there is a vertex $y_{r+2}^{r+2}$ such that $d_{Y_{r+2}}(\tilde Y_{r+2}^1\cup\tilde Y_{r+2}^2,y_{r+2}^{r+2})=d_{X_{r+1}}(\tilde X_{r+2}^1\cup\tilde X_{r+2}^2,x_{r+2}^{r+2})$ and the following properties hold. There exists an isomorphism $\tilde L_X\cup\tilde X_{r+2}^1\cup\tilde X_{r+2}^2\rightarrow \tilde L_Y\cup\tilde Y_{r+2}^1\cup\tilde Y_{r+2}^2$ which maps the vertices $x_{r+2}^1,\ldots,x_{r+2}^{r+2}$ to the vertices $y_{r+2}^1,\ldots,y_{r+2}^{r+2}$ respectively and the pair $(\tilde L_Y\cup\tilde Y_{r+2}^1\cup\tilde Y_{r+2}^2,\tilde Y_{r+2}^1\cup\tilde Y_{r+2}^2)$ is cyclically $2^{k-r-2}$-maximal, where $\tilde L_Y$ is a minimal path which connects $y_{r+2}^{r+2}$ and $\tilde Y_{r+2}^1\cup\tilde Y_{r+2}^2$ in $Y_{r+2}$. If $r=k-2$, then Duplicator wins. If $r<k-2$, then, obviously, the graphs $\tilde X_{r+2}^1:=\tilde X_{r+2}^1\cup\tilde X_{r+2}^2\cup\tilde L_X$ and $\tilde Y_{r+2}^1:=\tilde Y_{r+2}^1\cup\tilde Y_{r+2}^2\cup\tilde L_Y$ are $(k,r+2,1)$-regular equivalent in $(X_{r+2},Y_{r+2})$. Therefore, in the $r+3$-th round Duplicator exploits the strategy {\bf\sf SF}. Finally, if the vertex $x_{r+2}^{r+2}$ is not in the cyclic $2^{k-r-1}$-extension of the graph $\tilde X_{r+2}^2$ and $d_{X_{r+2}}(x_{r+2}^{r+2},\tilde X_{r+2}^1\cup\tilde X_{r+2}^2)> 2^{k-r-2}$, then denote the only (if it exists) cyclic $2^{k-r-2}$-extension of the graph $(\{x_{r+2}^{r+2}\},\varnothing)$ by $\tilde X_{r+2}^3$ (if there are no such extensions, then set $\tilde X_{r+2}^3=(\{x_{r+2}^{r+2}\},\varnothing)$). Obviously, $d_{X_{r+2}}(\tilde X_{r+2}^3,\tilde X_{r+2}^1\cup\tilde X_{r+2}^2)>2^{k-r-2}$. By Property {\sf 2)}, in $Y_{r+2}$ there are a vertex $y_{r+2}^{r+2}$ and a subgraph $\tilde Y_{r+2}^3$ such that $d_{Y_{r+2}}(\tilde Y_{r+2}^3,\tilde Y_{r+2}^1\cup\tilde Y_{r+2}^2)=2^{k-r-2}+1$, there exists an isomorphism $\tilde X_{r+2}^3\rightarrow\tilde Y_{r+2}^3$ which maps the vertex $x_{r+2}^{r+2}$ to the vertex $y_{r+2}^{r+2}$, and there are no cyclic $2^{k-r-2}$-extensions of the graph $\tilde Y_{r+2}^3$. If $r=k-2$, then Duplicator wins. If $r<k-2$, then, obviously, the ordered tuples $\tilde X_{r+2}^1, \tilde X_{r+2}^2,\tilde X_{r+2}^3$ and $\tilde Y_{r+2}^1,\tilde Y_{r+2}^2,Y_{r+2}^3$ are $(k,r+2,3)$-regular equivalent in $(X_{r+2},Y_{r+2})$. Therefore, in the $r+3$-th round Duplicator exploits the strategy {\bf\sf SF}.

\section{Extended law}

Theorem~\ref{new_th_3} can be extended in the following way.

\begin{theorem}
Let $k>3$, $b$ be arbitrary natural numbers. Moreover, let $\frac{a}{b}$ be an irreducible positive fraction, $\alpha=1-\frac{1}{2^{k-1}+a/b}$. Denote $\nu=\max\{1,2^{k-1}-b\}$. If $a\in\{\nu,\nu+1,\ldots,2^{k-1}\}$, then $\alpha\notin S_k^2$.
\label{spectrum}
\end{theorem}

A proof of the theorem is nearly the same as the proof of Theorem~8 from~\cite{Joint}, therefore, we do not give here a detailed proof. The idea is the following. As Duplicator has a winning strategy in the game EHR$(G,H,k)$ for all pairs of graphs $(G,H)$ such that $G,H\in\mathcal{S}$ (see Section~\ref{properties}), then by Theorem~\ref{ehren} it is sufficient to prove that for any $\alpha$ from the statement of Theorem~\ref{spectrum} there exists $\varepsilon$ such that ${\sf P}(G(n,p)\in\mathcal{S})\rightarrow 1$ as $n\rightarrow\infty$ for any $p\in[n^{-\alpha-\varepsilon},n^{-\alpha+\varepsilon}]$ (see the proof of Theorem~ 8 from~\cite{Joint}).

\section{Acknowledgements}

This work was carried out with the support of the Russian
Foundation for Basic Research grant No.~15-01-03530 and by the
grant No.~16-31-60052.

\end{document}